\documentclass[a4paper, 11pt]{article}
\usepackage{amsmath}
\usepackage{amsbsy}
\usepackage{amsfonts}
\usepackage{amstext}
\usepackage{amssymb}
\usepackage{MnSymbol}
\usepackage{mathbbol}
\usepackage{wasysym}
\usepackage{eucal}
\usepackage{mathtools}
\usepackage[amsmath,thmmarks]{ntheorem}		
\usepackage{stackrel}
\usepackage{ulem}
\usepackage{youngtab}
\usepackage{bold-extra}

\usepackage{times,tikz}
   \usetikzlibrary{matrix}
   \usetikzlibrary{calc}
   \usetikzlibrary{decorations}

\usepackage{float}
\usepackage{subfigure}

\usepackage{enumitem}

\usepackage[a4paper, left=2.1cm, right=2.1cm, top=3.0cm, bottom=3.5cm, bindingoffset=0cm]{geometry}
\usepackage{color}

\newcommand{\ZZ}{\ensuremath{{\mathbb Z}}}

\newcommand{\End}{\ensuremath{\text{End}}}

\newcommand{\spann}[1]{\ensuremath{\text{span}_{#1} \:}}
\renewcommand{\dim}{\ensuremath{\text{dim}}}

\newcommand{\ua}{\ensuremath{^{\ast}}}

\newcommand{\p}{\ensuremath{^{\prime}}}
\newcommand{\pp}{\ensuremath{^{\prime\prime}}}

\renewcommand{\t}[1]{\ensuremath{\tilde{#1}}}

\newcommand{\ha}[1]{\ensuremath{\hat{#1}}} 
\newcommand{\h}[1]{\ensuremath{\widehat{#1}}}

\newcommand{\ul}[1]{\ensuremath{\underline{#1}}}

\newcommand{\rar}{\rightarrow}

\newcommand{\zent}{\ensuremath{{\mathsf C}}}						

\newcommand{\antl}[1]{\ensuremath{\text{n}\h{\text{TL}}_{#1}}} 
\newcommand{\fntl}[1]{\ensuremath{\text{n}\text{TL}_{#1}}}		 
\newcommand{\atl}[1]{\ensuremath{\h{\text{TL}}_{#1}}}					 
\newcommand{\V}{\ensuremath{\mathsf{V}}} 
\newcommand{\W}{\ensuremath{\mathsf{W}}}

\newcommand{\ttt}{\ensuremath{\mathbf{t}}}

\newcommand{\ty}[1]{\ensuremath{\mathsf{#1}}} 

\newcommand{\seq}[1]{\ensuremath{\mathbf{\textsc{\Large#1}}}}

\newcommand{\vwi}{\ensuremath{v(\seq{i})}}
\newcommand{\vwj}{\ensuremath{v(\seq{j})}}
\newcommand{\vwip}{\ensuremath{v(\seq{i}\p)}} 
\newcommand{\ahj}{\ensuremath{a(\ha{\seq{j}})}}
\newcommand{\ahi}{\ensuremath{a(\ha{\seq{i}})}} 
\newcommand{\inj}{\ensuremath{\seq{i}^\text{in}_{\ul{j}}}}
\newcommand{\outj}{\ensuremath{\seq{i}^\text{out}_{\ul{j}}}}
\newcommand{\ellj}{\ensuremath{\ell_{\ul{j}}}} 
\newcommand{\inz}{\ensuremath{(\seq{i}^\text{in})^\ZZ}}
\newcommand{\outz}{\ensuremath{(\seq{i}^\text{out})^\ZZ}}
\newcommand{\vwedge}{\ensuremath{\bigoplus\limits_{k=0}^N\left(\ground[q] \otimes \bigwedge^k\ground^N\right)}}
\newcommand{\subsci} {\ensuremath{\mathbf{\textsc{I}}}}
\newcommand{\subscj} {\ensuremath{\mathbf{\textsc{J}}}}
\newcommand{\cl}{\ensuremath{\xi}}

\newcommand{\ground}{\ensuremath{\mathbb k}} 

\theoremstyle{plain}
\theoremseparator  {.}           
\theoremheaderfont {\bfseries}
\theorembodyfont   {\normalfont}
\theoremnumbering  {arabic}
\makeatletter
\newtheoremstyle{normal}%
{\item[\hskip\labelsep \theorem@headerfont ##1\ ##2\theorem@separator]\normalfont}%
{\item[\hskip\labelsep \theorem@headerfont ##1\ ##2]{\theorem@headerfont (##3)}\theorem@separator\ \normalfont}
\newtheoremstyle{nonumber}%
{\item[\theorem@headerfont\hskip\labelsep ##1\theorem@separator]\normalfont}%
{\item[\theorem@headerfont\hskip \labelsep ##1]{\theorem@headerfont (##3)}\theorem@separator\ \normalfont}
\makeatother

\theoremstyle{normal}
\newtheorem{thm}      {Theorem} [section]		

\newtheorem{lemma}      [thm]  {Lemma}        
 
\newtheorem{cor}        [thm]  {Corollary}

\newtheorem{prop}       [thm]  {Proposition}
\newtheorem{bem}        [thm]  {Remark}
\theoremsymbol{}
\newtheorem{bsp}        [thm]  {Example}

\newtheorem{defi}       [thm]  {Definition}

\theoremstyle{nonumber}
\theoremsymbol{$\square$} 
\newtheorem{bew}{Proof}
\theoremsymbol{}    

\theoremstyle{nonumber}
\newtheorem{thmnn}{Theorem}

\theoremstyle{nonumber}
\newtheorem{propnn}{Proposition}

\theoremstyle{nonumber}
\theoremstyle{nonumber}
\theoremstyle{nonumber}
\theoremsymbol{$\star$}

\theoremsymbol{}

\begin{document}

\title{The center of the affine nilTemperley-Lieb algebra}
\date{\today}
\author{Georgia Benkart$^*$ and Joanna Meinel\thanks{The authors thank Catharina Stroppel
for many helpful discussions and express their gratitude to
the Mathematical Sciences Research Institute (MSRI) where their joint research began. The second author would like to thank Daniel Tubbenhauer for remarks on the embeddings. This work is part of the second author's PhD project at the MPIM Bonn, supported by a grant of the Deutsche Telekom Stiftung.}}
\maketitle

\begin{abstract}
We give  a description of  the center of the affine nilTemperley-Lieb algebra based on a certain grading
of the algebra and on a faithful representation of it on fermionic particle configurations.   We present a normal form for monomials, hence construct a basis of the algebra, and use this basis  to show that the affine nilTemperley-Lieb algebra is finitely generated over its center.  As an application,  we obtain a natural embedding of the affine nilTemperley-Lieb algebra on $N$ generators into the affine nilTemperley-Lieb algebra on $N+1$ generators. 
\end{abstract}

\section{Introduction}

The main goal of this work is to describe the center of the affine nilTemperley-Lieb algebra $\antl{N}$ over any ground field. Only two tools are used: \ a fine grading on $\antl{N}$ and a representation of $\antl{N}$ on fermionic particle configurations on a circle.  It is essential that this graphical representation is faithful (see \cite[Prop.~9.1]{ks}).  
We provide an alternative proof of that fact  by constructing a basis for $\antl{N}$ that is especially adapted to the problem. This basis has further advantages:   It can be used to prove that the affine nilTemperley-Lieb algebra is finitely generated over its center.  Also, it can be used to exhibit an explicit embedding of $\antl{N}$ into $\antl{N+1}$
defined on basis elements that otherwise would not be apparent,  since the defining relations of these algebras are affine, and there is no embedding of the corresponding Coxeter graphs.

For a ground field $\ground$, the {\it affine nilTemperley-Lieb algebra} $\antl{N}$ is the unital associative $\ground$-algebra given by $N$ generators $a_0,\ldots,a_{N-1}$ and nil relations $a_i^2=0$ and $a_ia_{i\pm1}a_i=0$ for all $i$. Generators that are far apart commute, i.e. $a_ia_j=a_ja_i$ for $i-j \neq \pm 1 \text{ mod }N$. In these relations, the indices
are interpreted modulo $N$ so that the generators $a_0$ and $a_{N-1}$ are neighbours that do not commute.
The subalgebra of $\antl{N}$ generated by $a_1,\ldots,a_{N-1}$ is the {\it (finite) nilTemperley-Lieb algebra} $\fntl{N}$. 
The affine nilTemperley-Lieb algebra appears in many different settings, which we describe next. 
\begin{itemize}

\item $\antl{N}$ {\it is a quotient of the affine nilCoxeter algebra of type} $\t{\ty{A}}_{N-1}$.

\hspace{.3cm} The affine nilCoxeter algebra $\h{\mathsf U}_N$ of type $\t{\ty{A}}_{N-1}$ over a field $\ground$
is the unital associative algebra generated by elements $u_i$, $0 \leq i \leq N-1$, satisfying the 
relations $u_i^2 = 0$; $u_i u_j = u_j u_i$ for $i - j \neq \pm 1\text{ mod }N$;  and $u_i u_{i+1}u_i = u_{i+1}u_i u_{i+1}$ for $1 \leq i \leq N-1$, where the subscripts are read modulo  $N$.  The algebra $\antl{N}$ is isomorphic to the quotient
of $\h{\mathsf U}_N$ obtained by imposing the additional relations $u_i u_{i+1}u_i = u_{i+1}u_i u_{i+1}=0$ for
$1 \leq i \leq N-1$.    The affine nilCoxeter algebra is closely connected with affine Schur functions, $k$-Schur functions, and the affine Stanley symmetric functions, which are related to reduced word decompositions in the affine symmetric group (see e.g. \cite{L1}, \cite{L2}). 

\hspace{.3cm} The  nilCoxeter algebra $\mathsf{U}_{N}$ has generators $u_i, 1\leq i \leq N-1$, which satisfy the same  relations as they do in $\h{\mathsf U}_N$.  It  first appeared  in work on the cohomology of flag varieties  \cite{BGG} 
and  has played an essential role in studies on Schubert polynomials,
Stanley symmetric functions,  and the geometry of flag varieties (see for example \cite{LS}, \cite{M}, \cite{KK} \cite{FS}).   
The definition  of  $\mathsf{U}_N$  was inspired by the divided difference operators $\partial_i$ on polynomials in variables
$\mathbf{x} = \{x_1,  \dots,  x_{N}\}$  defined  by 
$$\partial_i(f) = \frac{f(\mathbf x) - f(\sigma_i \mathbf{x})}{x_i - x_{i+1}},$$  where the transposition
$\sigma_i$ fixes all the variables except for $x_i$ and $x_{i+1}$, which it interchanges.   
The operators $\partial_i$  satisfy the nilCoxeter relations above,  and applications of these relations enabled Fomin and Stanley \cite{FS}  to  recover known properties and establish new properties of Schubert polynomials.    

\hspace{.3cm}The algebra $\mathsf{U}_N$ belongs to a two-parameter family of algebras having
generators  $u_i$, $1 \leq i \leq N-1$,  which satisfy the relations $u_i u_j = u_j u_i$ for $\vert i-j \vert > 1$  and $u_i u_{i+1}u_i = u_{i+1}u_i u_{i+1}$ for $1 \leq i \leq N-2$ from above,  together with the relation 
 $u_i^2= \alpha u_i + \beta$ for all $i$, where $\alpha, \beta$ are fixed parameters.     
In particular, the specialization $\alpha=\beta=0$ yields the nilCoxeter algebra;
 $\alpha=0$, $\beta=1$ gives the standard presentation of the group algebra of the symmetric group $\ground \mathsf{S}_N$;  and  $\alpha= 1-q, \beta = q$ gives the Hecke algebra $\mathsf{H}_N(q)$ of type $\ty{A}$.

\hspace{.3cm} Motivated by categorification results in \cite{CF}, Khovanov \cite{Kh} introduced restriction and induction functors  $\mathsf{F}_D$ and $\mathsf{F}_{X}$ corresponding to the natural inclusion of 
algebras $\mathsf{U}_N \hookrightarrow \mathsf{U}_{N+1}$
on the direct sum $\mathcal C$ of the categories $\mathcal C_N$ of finite-dimensional $\mathsf{U}_N$-modules.  These functors categorify the Weyl  algebra of differential operators 
with polynomial coefficients in one variable and correspond to the Weyl algebra generators $\partial$ and $x$ (derivative and multiplication by $x$),  which satisfy the relation  $\partial x - x\partial  = 1$.  

\hspace{.3cm} Brichard \cite{bri} used a diagram calculus on cylinders  to determine the dimension of the center of $\mathsf{U}_N$ and to describe a basis of the center  for which the multiplication is trivial. In this diagram calculus on $N$ strands, the generator $u_i$ corresponds to a crossing of the strands $i$ and $i+1$.
The nil relation $u_i^2=0$ is represented by demanding that any two strands may cross at most once;  otherwise the diagram is identified with zero. 

\item $\antl{N}$ {\it is a quotient of the negative part of the universal enveloping algebra of the affine Lie algebra $\h{\mathfrak {sl}}_N$.}

\hspace{.3cm}  The negative part $U^-$  of the universal enveloping algebra $U$ 
of the affine Lie algebra $\h{\mathfrak{sl}}_N$  has generators $f_i$,  $0 \leq i \leq N-1$, which satisfy the Serre relations
$f_i^2f_{i+ 1} -  2 f_i f_{i+ 1}f_i + f_{i + 1}f_i^2 = 0 = f_{i+1}^2f_{i} -  2 f_{i+ 1} f_{i}f_{i+ 1}+ f_{i }f_{i+ 1}^2$
and $f_i f_j = f_j f_i$ for $i-j \neq \pm 1\text{ mod }N$ (all indices modulo $N$).    Factoring $U^-$ by the ideal generated by the elements $f_i^2$, $0 \leq i \leq N-1$,
gives $\antl{N}$ whenever the characteristic of $\ground$ is different from 2.  

\item $\antl{N}$  {\it acts on the small quantum cohomology ring of the Grassmannian.}   

\hspace{.3cm}  As in \cite[Sec.~2]{post}, (see also \cite{ks}),  consider the cohomology ring $\text{H}^\bullet(\text{Gr}(k,N))$ with integer coefficients  for the Grassmannian $\text{Gr}(k,N)$ of $k$-dimensional subspaces of $\ground^N$.  It has a basis given by the Schubert classes $[\Omega_\lambda]$,  where $\lambda$ runs over all partitions  with $k$ parts, 
the largest part having size $N-k$.  By recording the $k$ vertical and $N-k$ horizontal steps that identify the Young diagram of $\lambda$ inside the northwest corner of a $k \times (N-k)$ rectangle, such a partition corresponds to a $(0,1)$-sequence of length $N$ with $k$ ones (resp. $N-k$ zeros) in the positions corresponding to the 
vertical (resp. horizontal) steps.

\hspace{.3cm}  As a $\ZZ[q]$-module for an indeterminate $q$, the quantum cohomology ring of the Grassmannian is given by
$\text{qH}^\bullet(\text{Gr}(k,N))\ =\ \ZZ[q]\otimes_\ZZ\text{H}^\bullet(\text{Gr}(k,N))$
together with a $q$-multi\-plication.  The $\antl{N}$-action can be defined combinatorially on
$$\text{qH}^\bullet(\text{Gr}(k,N))\ \cong\ \spann{\ZZ[q]}\{(0,1)\text{-sequences of length }N\text{ with }k\text{ ones}\}
$$
as described in the next item,  and the multiplication of two Schubert classes $[\Omega_\lambda]\cdot[\Omega_\mu]$ is equal to $s_\lambda\cdot[\Omega_\mu]$ where $s_\lambda$ is a certain Schur polynomial in the generators of $\antl{N}$ as  in \cite[Cor.~8.3]{post}.

\item $\antl{N}$ \ {\it acts faithfully on fermionic particle configurations on a circle.} 

\hspace{.3cm} This is the graphical representation from \cite{ks} (see also \cite{post}), which we use  in our description of the center of $\antl{N}$.
First,  a $(0,1)$-sequence with $k$ ones is identified with a circular particle configuration having $N$ positions, where the $k$ particles are distributed at the position on the circle that corresponds to their position in the sequence,  so that there is at most one particle at each position.  On the space
$$\spann{\ground[q]}\{\text{fermionic particle configurations of }k\text{ particles on a circle with }N\text{ positions}\},$$
the generators $a_i$ of $\antl{N}$ act by sending a particle lying at position $i$ to position $i+1$. Additionally, the particle configuration is multiplied by $\pm q$ when applying $a_0$. The precise definition is given in Section~\ref{sec:rep}, but here is a representative picture:
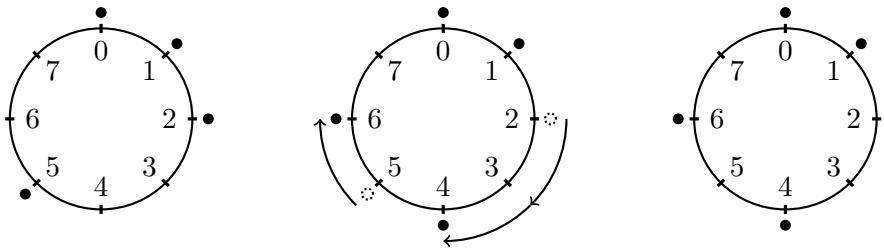
\begin{figure}[H]\label{fig:intro}
\begin{center}
\begin{tikzpicture}[scale=0.6]

\begin{scope}[xshift=-7.5cm]
\path[use as bounding box] (-3.5,-1.5) rectangle (3.5,2.5);
\draw[thick] (0,0) circle (2cm);

\foreach \x in {45,90,...,360}
\draw[very thick] (\x:1.9cm) -- (\x:2.1cm);

\foreach \x in {0,...,7}
\draw ({90-\x*45}:1.5cm) node {$\x$};

\foreach \x in {0,1,2,5}
\fill ({90-\x*45}:2.35cm) circle (1.2mm);
\end{scope}

\begin{scope}[xshift=0cm]
\path[use as bounding box] (-3.5,-1.5) rectangle (3.5,2.5);
\draw[thick] (0,0) circle (2cm);

\foreach \x in {45,90,...,360}
\draw[very thick] (\x:1.9cm) -- (\x:2.1cm);

\foreach \x in {0,...,7}
\draw ({90-\x*45}:1.5cm) node {$\x$};

\foreach \x in {0,1,4,6}
\fill ({90-\x*45}:2.35cm) circle (1.2mm);

\foreach \x in {2,5}
\draw[thick, densely dotted] ({90-\x*45}:2.35cm) circle (1.2mm);

\foreach \x in {2,3,5}
\draw[thick, ->] ({90-\x*45}:2.7cm) arc ({90-\x*45}:45-\x*45:2.7cm);
\end{scope}

\begin{scope}[xshift=7.5cm]
\path[use as bounding box] (-3.5,-1.5) rectangle (3.5,2.5);
\draw[thick] (0,0) circle (2cm);

\foreach \x in {45,90,...,360}
\draw[very thick] (\x:1.9cm) -- (\x:2.1cm);

\foreach \x in {0,...,7}
\draw ({90-\x*45}:1.5cm) node {$\x$};

\foreach \x in {0,1,4,6}
\fill ({90-\x*45}:2.35cm) circle (1.2mm);
\end{scope}
\end{tikzpicture}
\end{center}
\caption{$N=8$: Application of $a_3a_2a_5$ to the particle configuration $(0,1,2,5)$ gives $(0,1,4,6)$.}
\end{figure}

\item $\antl{N}$ {\it appears as a subalgebra of the annihilation/creation algebra.}

\hspace{.3cm} The finite nilTemperley-Lieb algebra is a subalgebra of the Clifford algebra having generators \newline \hfil $\{\cl_i,\cl_i\ua\ |\ 0\leq i\leq N-1\}$ and relations $\cl_i\cl_j+\cl_j\cl_i=0$, \ $\cl_i\ua\cl_j\ua+\cl_j\ua\cl_i\ua=0$, \  $\cl_i\cl_j\ua+\cl_j\ua\cl_i=\delta_{ij}$. The Clifford generators $\cl_i$ (resp. $\cl_i\ua$) act on the fermionic particle configurations by annihilation (resp. creation) of a particle at position $i$. The finite nilTemperley-Lieb algebra appears inside the Clifford algebra via $a_i\mapsto \cl_{i+1}\ua\cl_i$.  As discussed in \cite[Sec.~8]{ks},  the affine nilTemperley-Lieb algebra is a $q$-deformation of this construction. 

\item $\antl{N}$ {\it is the associated graded algebra of the affine Temperley-Lieb algebra.} 

\hspace{.3cm}  The affine Temperley-Lieb algebra $\atl{N}(\delta)$ has the usual commuting relations and the relations $a_ia_{i\pm1}a_i=a_i$ and  $a_i^2=\delta a_i$ for some parameter $\delta\in\ground$ instead of the nil relations
(where again all indices are mod $N$). It is a filtered algebra with its $\ell$th filtration space generated by all monomials of length $\leq \ell$.  Since its associated graded algebra is $\antl{N}$ for any value of $\delta$, elements of $\antl{N}$ can be identified with reduced expressions in $\atl{N}(\delta)$.

\hspace{.3cm}  The diagram algebra structure of $\atl{N}(\delta)$ is given by the same pictures as for the Temperley-Lieb algebra, but now the diagrams are wrapped around the cylinder (see e.g. \cite{fg}, \cite{kx}). The top and bottom of the cylinder each have $N$ nodes.  Monomials in the affine Temperley-Lieb algebra are represented by diagrams of $N$ non-crossing strands, each connecting a pair of those $2N$ nodes. Multiplication of two monomials is realized by stacking the cylinders one on top of the other, and connecting and smoothing the strands. Whenever  the strands form a circle, this is removed from the diagram at the expense of multiplying by  the parameter $\delta$. The relation $a_ia_{i\pm1}a_i=a_i$ corresponds to the isotopy between a strand that changes direction and a strand that is pulled straight.

\hspace{.3cm}  In contrast, the affine \emph{nil}Temperley-Lieb algebra is \emph{not} a diagram algebra.   The relation $a_ia_{i\pm1}a_i=0$ implies that isotopy would identify zero and nonzero elements. Nevertheless, the diagram of a reduced expression in $\atl{N}$ may be considered as an element of $\antl{N}$.  Such a diagram consists of a number (possibly 0)  of arcs that connect two nodes on the top of the cylinder, the same number of arcs connecting two nodes on the bottom, and arcs that connect a top node and a bottom one.  The latter arcs wrap around the cylinder either all in a strictly clockwise direction or all in a strictly counterclockwise way. Since the multiplication of two such diagrams may give zero, we will not use this diagrammatic realization here.
\end{itemize}

We proceed as follows:
In Section \ref{sec:notation},  we introduce the notation used in this article. 
The $\ZZ^N$-grading of $\antl{N}$ is given is Section \ref{sec:gradings},  and its importance for the description of the center is discussed. 
In Section \ref{sec:rep},  we give a detailed definition of the $\antl{N}$-action on particle configurations on a circle. We also define special monomials that serve as the projections onto a single particle configuration (up to multiplication by $\pm q$).  Proposition \ref{prop:faith} of that section recalls
\cite[Prop.~9.1]{ks} stating that the representation is faithful.  In \cite{ks}, this fact is deduced from the finite nilTemperley-Lieb algebra case,  as treated in \cite{bjs} and  \cite[Prop.~2.4.1]{bfz}. We give a complete, self-contained proof in Section \ref{sec:appendix}. Our proof is elementary and relies on the construction of a basis. 
In Section \ref{sec:center},  we state the main result (Theorem \ref{thm:center})  of this article:
\begin{thmnn}
The center of $\antl{N}$ is the subalgebra
$$\mathsf{C}_N = \mathsf{Cent}(\antl{N})\ =\ \langle 1, \ttt_1,\ldots, \ttt_{N-1}\rangle\ \cong\ \frac{\ground[\ttt_1,\ldots,\ttt_{N-1}]}{(\ttt_k\ttt_\ell \mid  k \neq \ell)},$$
where the generator $\ttt_k = (-1)^{k-1}\sum_{\mid \subsci \mid = k} \limits a(\ha{\seq{i}})$ is the sum of monomials $a(\ha{\seq{i}})$ corresponding to particle configurations given by increasing sequences $\seq{i} = \{1 \leq i_1 < \ldots < i_k \leq N\}$ of length $k$.  The  monomial $a(\ha{\seq{i}})$  sends particle configurations with $n \neq k$ particles to 0 and acts on a particle configuration with $k$ particles by  projecting onto $\seq{i}$ and multiplying by $(-1)^{k-1}q$.  Hence,
$\ttt_k$ acts as multiplication by $q$ on the configurations with $k$ particles. 
\end{thmnn}

Our $N-1$ central generators $\ttt_k$  are essentially the $N-1$ central elements constructed by Postnikov.  Lemma 9.4 of \cite{post} gives an alternative description of $\ttt_k$ as product of the $k$th elementary symmetric polynomial (with factors cyclically ordered) with the $(N-k)$th complete homogeneous symmetric polynomial (with factors reverse cyclically ordered) in the noncommuting generators of $\antl{N}$.   The above theorem shows that  in fact  these elements generate the entire center of $\antl{N}$.
In Section \ref{sec:finite},  we establish that $\antl{N}$ is finitely generated over its center.
In Section \ref{sec:inclusion},  we define a monomial basis for $\antl{N}$  indexed by pairs of particle configurations together with a natural number indicating how often the particles have been moved around the circle.  A proof that this is indeed a basis of $\antl{N}$ can be found in Section \ref{sec:appendix}. With this basis at hand,  we obtain inclusions $\antl{N}\subset\antl{N+1}$. The inclusions are not as obvious as those for the nilCoxeter algebra $\mathsf{U}_N$ having underlying Coxeter graph of type $\ty{A}_{N-1}$,  since one cannot deduce them from embeddings of the affine Coxeter graphs. Our result,  Theorem \ref{thm:emb},  reads as follows:
\begin{thmnn}
For all $0\leq m \leq N-1$, there are unital algebra embeddings $\varepsilon_m :\ \antl{N}\rar\antl{N+1}$ given by 
\begin{align*} 
a_i\ \mapsto\ a_i\  \ \text{ for}\ \, 0\leq i\leq m-1,\qquad a_{m}\ \mapsto\ a_{m+1} a_{m},\qquad a_i\ &\mapsto\ a_{i+1}\ \  \text{ for}\ \,  m+1\leq i\leq N-1.
\end{align*}\end{thmnn}
In Section \ref{sec:appendix},  we show how to construct the monomial basis, namely by using a normal form algorithm that reorders the factors of a nonzero monomial. Our basis is reminiscent of the Jones normal form for reduced expressions of monomials in the Temperley-Lieb algebra,  as discussed in \cite{rsa},  and is characterised in Theorem \ref{thm:basis} as follows: (See also Theorem \ref{thm:easybasis} which gives a different description.)  
\begin{thmnn}[Normal form]
Every nonzero monomial in the generators $a_j$ of $\antl{N}$ can be rewritten uniquely in the form 
$$ (a^{(m)}_{i_1}\ldots a^{(m)}_{i_k})\ldots(a^{(n+1)}_{i_1}\ldots a^{(n+1)}_{i_k})(a^{(n)}_{i_1}\ldots a^{(n)}_{i_k})\ldots(a^{(1)}_{i_1}\ldots a^{(1)}_{i_k})(a_{i_1}\ldots a_{i_k})$$
with $a^{(n)}_{i_\ell}\in\{1,a_0,a_1,\ldots,a_{N-1}\}\text{ for all }1\leq n\leq m,\ 1\leq \ell\leq k$,
such that
$$a^{(n+1)}_{i_\ell}\in\begin{cases}\{1\}\quad &\text{if }a^{(n)}_{i_\ell}=1,\\ \{1,a_{j+1}\}\quad&\text{if }a^{(n)}_{i_\ell}=a_j.\end{cases}$$ 
The factors $a_{i_1},\ldots,a_{i_k}$ are determined by the property that the generator $a_{i_\ell-1}$ does not appear to the right of $a_{i_\ell}$ in the original presentation of the monomial.
Alternatively, every nonzero monomial is uniquely determined by the following data from its action on the graphical representation:
\begin{itemize}
\item  the input particle configuration with the minimal number of particles on which it acts nontrivially,
\item  the output particle configuration,
\item the power of $q$ by which it acts.
\end{itemize}
\end{thmnn}

For the proof of this result,  we recall a characterisation of the nonzero monomials in $\antl{N}$ from \cite{green}. Then we prove faithfulness of the graphical representation of $\antl{N}$ by describing explicitly the matrices representing our basis elements.   Al Harbat \cite{ah} has recently described a normal form for fully commutative elements of the affine Temperley-Lieb algebra, which gives a different normal form when passing to $\antl{N}$. 

Our results hold over an arbitrary ground field $\ground$, even one of characteristic $2$, simply by ignoring signs in that case. In fact, our arguments work for any associative unital ground ring $R$ by replacing $\ground$-vector spaces and $\ground$-algebras with  free $R$-modules and $R$-algebras, respectively. In particular, the affine nilTemperley-Lieb algebra over $\ground$ is replaced by the $R$-algebra with the same generators and relations, and the polynomial ring $\ground[q]$ is replaced by $R[q]$. The ring $R$
is not required to be a domain or be commutative.  This is possible because our arguments mainly rely on investigating monomials in the generators of $\antl{N}$. 
However, for simplicity we have chosen to assume $\ground$ is a field throughout the article.

\section{Notation}\label{sec:notation}
Let $\ground$ be any field,  and assume $N$ is a positive integer.
The {\it affine nilTemperley-Lieb algebra} $\antl{N}$ of rank $N$ is the unital associative $\ground$-algebra generated by elements $a_0,\ldots,a_{N-1}$ subject to the relations
\begin{align*}
a_i^2\ &=\ 0\quad &&\text{for all }\ 0\leq i\leq N-1,\\
a_ia_j\ &=\ a_ja_i\quad &&\text{for all } \  i-j \neq \pm 1 \text{ mod }N,\\
a_ia_{i+1}a_i\ &=\ a_{i+1}a_ia_{i+1}\ =\ 0\quad&&\text{for all }\ 0\leq i\leq N-1,
\end{align*}
where all indices are taken modulo $N$, so in particular $a_{N-1}a_0a_{N-1}=a_0a_{N-1}a_0=0$. The {\it finite nilTemperley-Lieb algebra} $\fntl{N}$ is the subalgebra of $\antl{N}$ generated by $a_1,\ldots,a_{N-1}$. 
We adopt the convention that $\fntl{1} =\ground 1$.  We fix the following notation for monomials in $\antl{N}$ and $\fntl{N}$: For an ordered index sequence $\ul{j}=(j_1,\ldots, j_m)$ with $0\leq j_1,\ldots,j_m\leq N-1$, we define the ordered monomial $a(\ul{j})= a_{j_1}\ldots a_{j_m}$.
Unless otherwise specified, we use the letters $i,j$ for indices from $\ZZ/N\ZZ$;  in particular, we often identify the indices $0$ and $N$.   

{\it Throughout we will assume $N\geq3$.}    

\section{Gradings}\label{sec:gradings}
One of the ingredients needed in Section \ref{sec:center} to study the center of $\antl{N}$ is a fine grading on the algebra.
Gradings faciliate the computation of the center of an algebra, as the following standard result reduces the work to determining homogeneous central elements.
\begin{lemma}\label{lem:tech1}
Let $A=\bigoplus\limits_{g\in G} A_g$ be an algebra graded by some abelian group $G$.
The center of $A$  is homogeneous, i.e. it inherits the grading.
\end{lemma}

\begin{bew}
Let $a=\sum\limits_{g\in G}a_g$ be a central element of  the graded algebra $A=\bigoplus\limits_{g\in G} A_g$.  We have for $b_h\in A_h$ that
$\sum\limits_{g\in G} a_g b_h = a b_h = b_h a= \sum\limits_{g\in G} b_h a_g.$
Since this equality must hold in every graded component, we get $a_g b_h = b_h a_g$ for all homogeneous elements $b_h$. Now take any element $b=\sum\limits_{h\in G}b_h$ in $A$, then
$a_g b = \sum\limits_{h\in G}a_g b_h = \sum\limits_{h\in G}b_h a_g= b a_g,$
hence $a_g$ is central.
\end{bew}

Since the defining relations are homogeneous, both $\antl{N}$ and $\fntl{N}$ have a $\ZZ$-grading by the length of a monomial, i.e. all generators $a_i$ have $\ZZ$-degree $1$.
This can be refined to a $\ZZ^N$-grading by assigning to the generator $a_i$ the degree $\zeta_i$, the $i$th standard basis vector in $\ZZ^N$. In either grading, we say that the degree $0$ part of an element in $\antl{N}$ or $\fntl{N}$ is its constant term.
\begin{bem}
Why do we exclude the case of $N\leq 2$ from our considerations?   For $N = 1,2$, there are isomorphisms $\antl{N}\cong \fntl{N+1}$,  and in these cases  the center is uninteresting.  The algebra $\antl{1}$  is $2$-dimensional and commutative; while $\antl{2}$ has dimension $5$,  and its center can be computed by hand making use of Lemma \ref{lem:tech1} and can be shown to be the $\ground$-span of $1, a_0a_1, a_1a_0$. 
\end{bem}
\begin{bem}
The affine (or finite) Temperley-Lieb algebra, which has relations $a_i a_j = a_j a_i$ for $i-j \neq \pm 1\ (\text{mod }N)$, $a_ia_{i\pm1}a_i=a_i$, 
and  $a_i^2=\delta a_i$ for some $\delta\in\ground$, is a filtered algebra with respect to the length filtration. For this algebra, the $\ell$th filtration space is generated by all monomials of length $\leq \ell$. Its associated graded algebra is $\antl{N}$ (or $\fntl{N}$). Thus, $\antl{N}$ is infinite dimensional when $N \geq 3$, while $\fntl{N}$ has dimension equal to the $N$th Catalan number $\frac{1}{N+1}{{2N} \choose N}$. \end{bem}

\section{A faithful representation}\label{sec:rep}
The second ingredient we use to determine the center is a faithful representation of $\antl{N}$. Here we recall the definition of the representation from \cite{ks} and describe its graphical realization, which is very convenient to work with.

Fix a basis $v_1,\ldots, v_N$ of $\ground^N$.  
Consider the vector space $\V=\vwedge$. It has a standard $\ground[q]$-basis consisting of wedges 
$$\vwi \ :=\ v_{i_1}\wedge\ldots\wedge v_{i_k}\quad\text{for all (strictly) increasing sequences } \ \seq{i}= \{1 \leq i_1 < \ldots <  i_k \leq N\}$$
for all $0\leq k\leq N$, where the basis element of $\ground=\bigwedge^0\ground^N$ is denoted $v(\emptyset)$. Throughout the rest of the article, all tensor products are taken over $\ground$, and we omit the tensor symbol in $\ground[q]$-linear combinations of wedges.
\begin{bem}
The indices of the vectors $v_j$ should be interpreted modulo $N$. We make no distinction between  $v_0$ and $v_N$
and often use the two interchangeably.
\end{bem}

It is helpful to visualize the basis elements $\vwi$ as particle configurations having $0\leq k\leq N$ particles arranged on a circle with $N$ positions, where there is at most one particle at each site, as pictured below for $N=8$ and $v(1,5,6)=v_1\wedge v_5\wedge v_6$. 
The vector $v(\emptyset)$ corresponds to the configuration with no particles. Then $\V$ is the $\ground[q]$-span of such circular particle configurations.
\begin{figure}[H]\label{fig:ex1}
\begin{center}
\begin{tikzpicture}[scale=0.6]
\draw[thick] (0,0) circle (2cm);

\foreach \x in {45,90,...,360}
\draw[very thick] (\x:1.9cm) -- (\x:2.1cm);

\foreach \x in {0,...,7}
\draw ({90-\x*45}:1.5cm) node {$\x$};

\foreach \x in {1,5,6}
\fill ({90-\x*45}:2.35cm) circle (1.2mm);

\end{tikzpicture}
\caption{The element $v_1\wedge v_5\wedge v_6$ in the graphical realization.}
\end{center}
\end{figure}
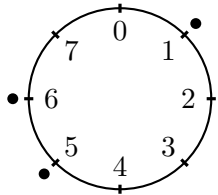

There is an action of the affine nilTemperley-Lieb algebra $\antl{N}$ defined on the basis vectors $\vwi$ of $\V$ as follows: 
\begin{defi}\label{defi:action}
For $1\leq j\leq N-1$,
\begin{align*}
a_j \vwi \ &=\ \begin{cases}
v_{i_1}\wedge\ldots\wedge v_{i_{\ell-1}}\wedge v_{j+1}\wedge v_{i_{\ell+1}} \wedge\ldots\wedge v_{i_k},\ &\text{if }i_\ell=j\text{ for some }\ell,\\
0,\quad &\text{otherwise.}\end{cases}
\end{align*}
For the action of $a_0$, note that $v_N$ appears in the basis element $\vwi$ if and only if it occurs in the last position, i.e. $v_{i_k}=v_N$, and define
\begin{align*}
a_0 \vwi \ &=\ \begin{cases}
q\cdot v_{i_1}\wedge\ldots\wedge v_{i_{k-1}}\wedge v_{1},\ &\text{if }i_k=N,\\
0,\quad &\text{otherwise;}\end{cases}\\
&=\ \begin{cases}
(-1)^{k-1}q\cdot v_1\wedge v_{i_1}\wedge\ldots\wedge v_{i_{k-1}},\ &\text{if }i_k=N,\\
0,\quad &\text{otherwise.}\end{cases}
\end{align*} \end{defi}
\begin{bem}
It follows that $a_j \vwi =0$ if the sequence $\seq{i}$ contains $j+1$ or if it does not contain $j$. 
In other words, $a_j$ acts by replacing $v_j$ by $v_{j+1}$. If this creates a wedge expression with two factors equal to $v_{j+1}$, the result is zero. In the graphical description, $a_j$ moves a particle clockwise from position $j$ to position $j+1$, and one records `passing position $0$' by multiplying by $\pm q$ as illustrated by the particle configurations
below. \end{bem}
\begin{figure}[H]\label{fig:ex2}
\begin{center}
\subfigure[$a_6(v_1\wedge v_5\wedge v_6)=v_1\wedge v_5\wedge v_7$]{\begin{tikzpicture}[scale=0.6]
\path[use as bounding box] (-3.5,-3) rectangle (3.5,3);
\draw[thick] (0,0) circle (2cm);

\foreach \x in {45,90,...,360}
\draw[very thick] (\x:1.9cm) -- (\x:2.1cm);

\foreach \x in {0,...,7}
\draw ({90-\x*45}:1.5cm) node {$\x$};

\foreach \x in {1,5,7}
\fill ({90-\x*45}:2.35cm) circle (1.2mm);

\draw[thick, densely dotted] (180:2.35cm) circle (1.2mm);
\draw[thick, ->] (-2.7cm,0cm) arc (180:135:2.7cm);

\end{tikzpicture}
}\qquad
\subfigure[$a_7a_1a_6(v_1\wedge v_5\wedge v_6)=v_2\wedge v_5\wedge v_0$]{\begin{tikzpicture}[scale=0.6]
\path[use as bounding box] (-3.9,-3) rectangle (3.9,3);
\draw[thick] (0,0) circle (2cm);

\foreach \x in {45,90,...,360}
\draw[very thick] (\x:1.9cm) -- (\x:2.1cm);

\foreach \x in {0,...,7}
\draw ({90-\x*45}:1.5cm) node {$\x$}; 

\foreach \x in {2,5,0}
\fill ({90-\x*45}:2.35cm) circle (1.2mm);

\draw[thick, densely dotted] (180:2.35cm) circle (1.2mm);
\draw[thick, ->] (180:2.7cm) arc (180:135:2.7cm);
\draw[thick, ->] (135:2.7cm) arc (135:90:2.7cm);
\draw[thick, densely dotted] (45:2.35cm) circle (1.2mm);
\draw[thick, ->] (45:2.7cm) arc (45:0:2.7cm);

\end{tikzpicture}
}\qquad
\subfigure[$a_0(v_5\wedge v_0) = -q\cdot v_1\wedge v_5$.]{\begin{tikzpicture}[scale=0.6]
\path[use as bounding box] (-3.5,-3) rectangle (3.5,3);

\draw[thick] (0,0) circle (2cm);

\foreach \x in {45,90,...,360}
\draw[very thick] (\x:1.9cm) -- (\x:2.1cm);

\foreach \x in {0,...,7}
\draw ({90-\x*45}:1.5cm) node {$\x$};

\foreach \x in {1,5}
\fill ({90-\x*45}:2.35cm) circle (1.2mm);

\draw[thick, densely dotted] (90:2.35cm) circle (1.2mm);
\draw[thick, ->] (90:2.7cm) arc (90:45:2.7cm);
\draw (70:3cm) node {$\cdot (-q)$};
\end{tikzpicture}
}
\caption{Examples for the action of $\antl{N}$ on a particle configuration}
\end{center}
\end{figure}
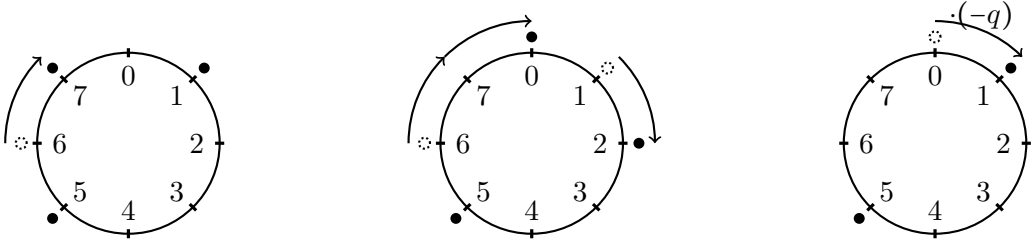

It is easy to verify that the defining relations for $\antl{N}$ hold for this action, assuming that $N\geq3$. Hence we obtain
\begin{lemma} \label{lem:constwedge}  \begin{itemize}
\item [{\rm (a)}]
Definition \ref{defi:action} gives a representation of $\antl{N}$ on $\V$.
\item [{\rm (b)}]
The number of wedges (i.e., the number of particles) remains constant under the action of the generators $a_i$, so that  $\V=\vwedge$ is a direct sum decomposition of $\V$ as an $\antl{N}$-module. \end{itemize}
\end{lemma}

The following crucial statement is taken from \cite[Prop.~2.4.1]{bfz} and \cite[Prop.~9.1.(2)]{ks}. We will give a detailed proof adapted to our notation in Section \ref{sec:appendix}.
\begin{prop}\label{prop:faith}
The action from Definition \ref{defi:action} gives a faithful representation of $\antl{N}$ on $\V$ when $N \geq 3$.
\end{prop}

From now on,  we will identify elements of $\antl{N}$ with their action on the particle configurations of the graphical representation.
\begin{bem}
The spaces $\ground[q]\otimes \bigwedge^0\ground^N$ and $\ground[q] \otimes \bigwedge^N\ground^N$ are trivial summands in $\V$ on which every generator $a_i$ acts as 0, and so they may be ignored 
when proving Proposition \ref{prop:faith}.
\end{bem}

For a standard basis element  $\vwi$  of  $1 \leq k \leq N-1$  wedges  corresponding to an increasing sequence $\seq{i} = \{1\leq i_1<\ldots< i_k\leq N\}$,  the next lemma defines a certain monomial $\ahi$ that projects  $\vwi$ onto $(-1)^{k-1} q\,\vwi$ and sends $\vwip$ to zero for $\seq{i}\p \neq \seq{i}$.
Before stating the result, we give an example to demonstrate in the graphical description how this projector will be defined.
\begin{bsp}\label{ex:part} 
Let $N=8$,  and consider the particle configuration $\vwi=v_1\wedge v_5\wedge v_6$.
With $a(\h{1\ 5\ 6})\ =\ (a_0 a_7)\cdot (a_4 a_3 a_2)\cdot (a_1 a_5 a_6)$
we obtain $a(\h{1\ 5\ 6})\cdot\ v_1\wedge v_5\wedge v_6\ =\ (-1)^2 q\cdot v_1\wedge v_5\wedge v_6$, which looks as follows in the graphical description:
\begin{figure}[H]
\begin{center}
\begin{tikzpicture}[scale=0.6]
\draw[thick] (0,0) circle (2cm);

\foreach \x in {45,90,...,360}
\draw[very thick] (\x:1.9cm) -- (\x:2.1cm);

\foreach \x in {0,...,7}
\draw ({90-\x*45}:1.5cm) node {$\x$};

\foreach \x in {6,5,1}
\fill[color=black] ({90-\x*45+10}:2.35cm) circle (1.2mm);
\foreach \x in {6,5,1}
\draw[thick, densely dotted, color=black] ({90-\x*45}:2.35cm) circle (1.2mm);

\draw[thick, ->] (177:2.7cm) arc (177:135:2.7cm);
\draw[thick, ->] (135:2.7cm) arc (135:90:2.7cm);
\draw[thick, ->] (90:2.7cm) arc (90:55:2.7cm);
\draw (70:3cm) node {$\cdot q$};

\draw[thick, ->] (42:2.7cm) arc (42:0:2.7cm);
\draw[thick, ->] (0:2.7cm) arc (0:-45:2.7cm);
\draw[thick, ->] (315:2.7cm) arc (315:270:2.7cm);
\draw[thick, ->] (270:2.7cm) arc (270:235:2.7cm);

\draw[thick, ->] (222:2.7cm) arc (222:190:2.7cm);

\end{tikzpicture}
\caption{The action of $a(\h{1\ 5\ 6})$ on the particle configuration $v_1\wedge v_5\wedge v_6$}
\end{center}
\end{figure}
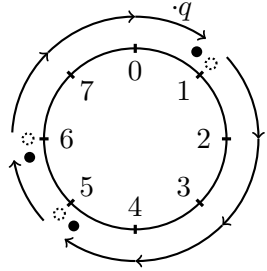  
The factor $a_1a_5a_6$ moves every particle one step forward clockwise. It is critical that we start by moving the particle at position $6$ before moving the particle at position $5$,  as otherwise the result would be zero. But since there is a `gap' at position $7$, we can move the particle from site $6$ to $7$, and afterwards the particle from site $5$ to $6$, without obtaining zero. The assumption that $k<N$ ensures such a gap always exists.

After applying $a_1a_5a_6$, the particles are at positions $2$, $6$,  and $7$.   The particle previously at position 5 is now at position 6, which is
where we want a particle to be.   The particle currently at position $2$  can be moved to position $5$ by applying  the product $a_4 a_3 a_2$. 
The particle now at position $7$ can be moved by $a_0 a_7$ to position $1$.  Hence, the result of applying $(a_0 a_7)\cdot (a_4 a_3 a_2) \cdot (a_1 a_5 a_6)$ is the same particle configuration as the original one. However, the answer must be multiplied by $\pm q$, since applying $a_0 a_7$  involves crossing the zero position once.   To determine the sign, note from Definition \ref{defi:action}  that  
$(a_0 a_7) \cdot (a_4 a_3 a_2) \cdot (a_1 a_5 a_6)(v_1 \wedge v_5 \wedge v_6) = q \cdot v_5 \wedge v_6 \wedge v_1
= (-1)^2 q \cdot v_1 \wedge v_5 \wedge v_6$, so the sign is $+$.  \end{bsp}

Now we describe the general procedure:
\begin{lemma}\label{lem:tech3}
Assume $\vwi$ is a particle configuration, where $\seq{i}= \{1\leq i_1<\ldots< i_k\leq N\}$ is an increasing sequence and $1 \leq k \leq N-1$.  
Then there exists an index $\ell$ such that $i_\ell+1<i_{\ell+1}$ (or  $i_{k}+1<i_1$), i.e. the sequence has a `gap' between $i_\ell$ and $i_{\ell+1}$.  Split the sequence $\seq{i}$ into the two parts $\{i_1<\ldots< i_\ell\}$ and $\{i_{\ell+1}<\ldots<i_k\}$.
Set
\begin{align*}
\ahi\ :=\ &(a_{i_{1}-1}a_{i_{1}-2}\ldots a_{i_k+2}a_{i_k+1}) \cdot \prod\limits_{s=1}^{k-1} (a_{i_{s+1}-1}a_{i_{s+1}-2}\ldots a_{i_s+2}a_{i_s+1})\tag{$\star$}\label{eq:ahidef}\\
& \qquad \cdot (a_{i_{\ell+1}}a_{i_{\ell+2}}\ldots a_{i_{k-1}}a_{i_{k}})\cdot (a_{i_1}a_{i_{2}}\ldots a_{i_{\ell-1}}a_{i_\ell}),
\end{align*}
where the indices are modulo $N$ in the factor  $(a_{i_{1}-1}a_{i_{1}-2}\ldots a_{i_k+2}a_{i_k+1})$. 
Then
$$\ahi \vwip = \begin{cases}(-1)^{k-1}q\cdot \vwi \quad &\text{if } \  \seq{i}\p=\seq{i},\\0\quad&\text{for all }\ \seq{i}\p\neq \seq{i} \ \text{ (of any length),}\end{cases}$$
and $\ahi$ has $\ZZ^N$-degree $(1,1,\ldots,1)$.
\end{lemma}
\begin{bew}
The assertions can be seen using the graphical realization of $\V$. The terms in the second line of equation \eqref{eq:ahidef} move a particle at site $i_j\in \seq{i}$ one step forward to $i_j+1$ for each $j$, while the terms in the first line send the particle from $i_j+1$ to the original position of  $i_{j+1}$. 

Consider first $\ahi \vwi$. By applying $(a_{i_{\ell+1}}a_{i_{\ell+2}}\ldots a_{i_{k-1}}a_{i_{k}})\cdot (a_{i_1}a_{i_{2}}\ldots a_{i_{\ell-1}}a_{i_\ell})$, every particle is first moved clockwise by one position.  By our choice of the index $i_\ell$,  we avoid mapping  the whole particle configuration  to zero. After that step, every particle is moved by one of the factors $(a_{i_{s+1}-1}a_{i_{s+1}-2}\ldots a_{i_s+2}a_{i_s+1})$ to the original position of its successor in the sequence $\seq{i}$, so the particle configuration remains the same. One of the particles has passed the zero position, so we have to multiply by $\pm q$.  Definition \ref{defi:action} tells us the appropriate sign is $(-1)^{k-1}$. 

Now consider $\ahi \vwip$ for $\seq{i}\p \neq \seq{i}$. The monomial $(a_{i_{\ell+1}}a_{i_{\ell+2}}\ldots a_{i_{k-1}}a_{i_{k}})\cdot (a_{i_1}a_{i_{2}}\ldots a_{i_{\ell-1}}a_{i_\ell})$ expects a particle at each of the sites $i_1,\ldots, i_k$, so if any of these positions is empty in $\vwip$, the result of applying $\ahi$ is zero. 
If the positions $i_1,\ldots, i_k$ are already filled, and there is an additional particle somewhere, multiplication by $(a_{i_{\ell+1}-1}a_{i_{\ell+1}-2}\ldots a_{i_\ell+2}a_{i_\ell+1})$ will cause two particles to be at the same position, hence the result is again zero.

Since every $a_j$ appears in $\ahi$ exactly once,  the monomial $\ahi$ has $\ZZ^N$-degree $(1,1, \ldots, 1)$.
\end{bew}
\begin{bsp} 
In the previous example, $N = 8$,  $\seq{i} = (1, 5, 6)$, and we may assume the two subsequences are $(1)$
and $(5,6)$.  Then the terms in the second line of 
\eqref{eq:ahidef} are  $(a_5 a_6)\cdot (a_1)  = a_1 a_5 a_6$.   The term corresponding to  $j = 1$ in the product on the first line of \eqref{eq:ahidef} 
is $a_4 a_3 a_2$,  and the expression corresponding to $j=2$ is empty, hence taken to be 1.  The first factor on the first line is $a_0 a_7$.  Thus, 
for $\seq{i}= (1, 5, 6)$,   $\ahi = (a_0 a_7) \cdot (a_4 a_3 a_2) \cdot (a_1 a_5 a_6)$, as in Example \ref{ex:part}. 
If the gap between $6$ and $0$ is used instead,  the right-hand factor of the second line is $a_1 a_5 a_6$ and
the left-hand factor is 1.  The factors in the first line remain the same, and so one obtains the same expression for $\ahi$.
\end{bsp}

\begin{bem}\label{bem:reconstruct}
Because $\V$ is a faithful module, $\ahi$ is, as an element in $\antl{N}$ (i.e. up to reordering according to the defining relations), uniquely determined by the increasing sequence $\seq{i}$.
One can read off $\seq{i}$ from $\ahi$ as follows: In the defining equation \eqref{eq:ahidef} of $\ahi$, the factors in the first line
 are pairwise commuting. The underlying subsequence $(i_{s+1}-1, i_{s+1}-2,\ldots, i_s+2, i_s+1)$  corresponding
 to the factor  $a_{i_{s+1}-1}a_{i_{s+1}-2}\ldots a_{i_s+2}a_{i_s+1}$ of $\ahi$ is a decreasing sequence. After all such decreasing sequences are removed from $\ahi$, what remains is a product of generators $a_j$ with an  increasing subsequence of
 indices or a product of two such subsequences  corresponding to the factors in the second line. This is $\seq{i}$. Given any monomial $a(\ul{r})$ of $\ZZ^N$-degree $(1,\ldots,1)$, one can rewrite it using the relations in $\antl{N}$ so that it is of the form $\ahi$ for some increasing sequence $\seq{i}$. Then $\vwi$ is the unique standard basis element upon which $a(\ul{r})=\ahi$ acts by multiplication by $\pm q$.
\end{bem}

\section{Description of the center}\label{sec:center}
In this section, we give an explicit description of the center $\zent_N$ of $\antl{N}$. We start with the following initial characterisation of the central elements:
\begin{lemma}\label{lem:tech2}
Any central element $c$ in $\antl{N}$ with constant term $0$ is a linear combination of monomials $a(\ul{j})=a_{j_1}\cdot\ldots \cdot a_{j_m}$ where every generator $a_i$, $0\leq i\leq N-1$, appears at least once. In particular,  a homogeneous nonconstant central element $c$ has $\ZZ$-degree at least $N$.
\end{lemma}
\begin{bew}
Assume $c = \sum\limits_{\ul{j}}c_{\ul{j}}a(\ul{j})$, where $c_{\ul{j}}\in\ground$ for all $\ul{j}$. By Lemma \ref{lem:tech1}, we can assume $c$ is a homogeneous central element with respect to the $\ZZ^N$-grading. By our assumption, $c\notin \ground$. For all $i$,  we need to show that $a_i$ occurs in each monomial $a(\ul{j})$ appearing in $c$.
Without loss of generality, we show this for $i=0$. 
Suppose some summand is missing $a_0$, then no summand contains $a_0$ because $c$ is homogeneous. Hence
$a_0 a(\ul{j}) \neq 0$ and $a(\ul{j})  a_0\neq0$ for all $\ul{j}$ with $c_{\ul{j}} \neq 0$, and since $a_0  c= c a_0$, none of the $a(\ul{j})$ can contain the factor $a_1$ either, as otherwise the factor $a_0$ cannot pass through $c$ from left to right (so also $a_{N-1}$ cannot be contained in the $a(\ul{j})$). Proceeding inductively, we see that all $a(\ul{j})$ must be a constant, contrary to our assumption.
\end{bew}

The next proposition states that on the standard wedge basis vector $\vwi$ of $\V$, any central element acts via multiplication by a polynomial $p_k\in \ground[q]$ that only depends on the length $k = \vert \seq{i}\vert$ of the 
increasing sequence  $\seq{i}=\{1 \leq i_1<\ldots < i_k\leq N\}$.
In other words, the decomposition of $\V$ into the summands $\ground[q] \otimes \bigwedge^k\ground^N$ is a decomposition with respect to different central characters (apart from the two trivial summands for $k\in\{0,N\}$).
\begin{prop}\label{prop:central}
For any central element $c\in\antl{N}$ and all increasing sequences $\seq{i}$ with fixed length $k$, there is some element $p_k\in\ground[q]$ such that $c \vwi =p_k\,\vwi$.
\end{prop}
\begin{bew}
We may assume $c$ is a nonconstant $\ZZ^N$-homogeneous central element of $\antl{N}$.
For $k\in\{0,N\}$, the action of  a generator  $a_i$ on a monomial of length $k$ is $0$, so $p_k=0$ for such values of $k$. Now consider $1\leq k\leq N-1$,  and suppose that $\seq{i}= \{1\leq i_1<\ldots< i_k\leq N\}$ is an increasing sequence of length $k$.
According to Lemma \ref{lem:constwedge}\,(b), the number of wedges in a vector remains constant under the action of the $a_i$. Hence 
$c \vwi =\sum\limits_{\mid\mathbf{\textsc{I}}\p\mid =k} c_{\mathbf{\textsc{I}} \p}\, \vwip$ for some polynomials $c_{\mathbf{\textsc{I}}\p}\in\ground[q]$.
We want  to prove that $c_{\mathbf{\textsc{I}}\p}=0$ for all $\seq{i}\p \neq \seq{i}$.

We have shown in Lemma \ref{lem:tech3} that to each increasing sequence $\seq{j}  \subset \{1,  \dots, N\}$ there corresponds a monomial $\ahj  \in\antl{N}$ that allows us to select a single basis vector:
$$\ahj \vwi \ =\ \begin{cases}(-1)^{k-1}q v(\seq{j}) \ &\text{if } \  \seq{i} =\seq{j},\\ 0\quad& \ \text{otherwise}.\end{cases}$$
Thus, for $\seq{j} \neq \seq{i}$,  we see that 
$$0\ =\ c (\ahj \vwi)\ =\ \ahj(c \vwi)\ =\ \ahj \left(\sum\limits_{\vert \mathbf{\textsc{I}}\p\vert=k} c_{\mathbf{\textsc{I}}\p}\,\vwip\right) = c_{\mathbf{\textsc{J}}}\,(-1)^{k-1}q v(\seq{j}),$$  implying  $c_{\mathbf{\textsc{J}}}=0$ for $\seq{j} \neq \seq{i}$. Hence, we may assume for each increasing sequence $\seq{i}$ that  $c \vwi = p_{\mathbf{\textsc{I}}}\,\vwi$ for
some polynomial $p_{\mathbf{\textsc{I}}} \in \ground[q]$.  Now it is left to show that $p_{\mathbf{\textsc{I}}} = p_{\mathbf{\textsc{I}}\p}$ for all $\seq{i}\p$ with  $\vert \seq{i}\p \vert = \vert \seq{i} \vert =k$.  It is enough to verify this for $\seq{i}$, $\seq{i}\p$ which differ in exactly one entry, i.e. $i_s=i$, $i_s\p=i+1$, and $i_\ell=i_\ell\p$ for all $\ell\neq s$, for some $1\leq s\leq k$ and $i\in\ZZ/N\ZZ$.
If  $1\leq i\leq N-1$,  we have
$$p_{\mathbf{\textsc{I}}\p}\,v(\seq{i}\p)\ =\ c v(\seq{i}\p) \ =\ c(a_i v(\seq{i}))  \ =\ a_i (c v(\seq{i})) \ = \ a_i (p_{\mathbf{\textsc{I}}}\,v(\seq{i}) ) \ = \ p_{\mathbf{\textsc{I}}}\,v(\seq{i}\p),$$
and if $i=0$,  we get
$$(-1)^{k-1}q p_{\mathbf{\textsc{I}}\p}\, v(\seq{i}\p)\ =\ (-1)^{k-1}q c v(\seq{i}\p) \ =\ c(a_0 v(\seq{i}))  \ =\ a_0 (cv(\seq{i})) \ = \ a_0 (p_{\mathbf{\textsc{I}}}\,v(\seq{i})) \ = \ (-1)^{k-1}q p_{\mathbf{\textsc{I}}}\,v(\seq{i}\p).$$
Hence,  $p_{\mathbf{\textsc{I}}\p}=p_{\mathbf{\textsc{I}}}$, and this common polynomial is the desired polynomial $p_k$.
\end{bew}
\begin{cor}
Any central element in $\antl{N}$ with constant term $0$ acts on a standard basis vector  $\vwi \in \V$ as multiplication by an element of $q\ground[q]$.
\end{cor}
\begin{bew}
According to Lemma \ref{lem:tech2}, each summand of such a central element must contain the factor $a_0$, and $a_0$ acts on a wedge product  by $0$ or multiplication by $\pm q$.
\end{bew}

Now we are ready to introduce  nontrivial central elements in $\antl{N}$.   For each $1\leq k\leq N-1$, set
\begin{equation} \ttt_k\ :=\ (-1)^{k-1}\sum\limits_{\mid \subsci \mid =k}\ahi,\end{equation}
where the monomials $\ahi$ correspond to increasing sequences $\seq{i} = \{1 \leq i_1 < \ldots < i_k \leq N\}$
of length $k$ as  defined in Lemma \ref{lem:tech3}.
\begin{bsp}
In $\antl{3}$:
\begin{align*}
\ttt_1\ &=\   a_2 a_1 a_0 + a_0 a_2 a_1 + a_1 a_0 a_2,\\
\ttt_2\ &=\ -a_0a_1a_2 - a_1a_2a_0 - a_2 a_0 a_1.
\end{align*}
In $\antl{4}$:
\begin{align*}
\ttt_1\ &=\ a_3a_2 a_1 a_0 + a_0 a_3a_2 a_1 + a_1 a_0a_3 a_2 + a_2 a_1 a_0 a_3,\\
\ttt_2\ &=\ -a_0 a_2 a_1 a_3 - a_1 a_3 a_0 a_2 - a_0 a_1 a_3 a_2 - a_1 a_2 a_0 a_3 - a_2 a_3 a_1 a_0 - a_3 a_0 a_2 a_1,\\
\ttt_3\ &=\ a_0a_1a_2a_3 + a_1a_2a_3a_0 + a_2a_3 a_0 a_1 + a_3a_0a_1 a_2.
\end{align*}
\end{bsp}

In the graphical realization of $\V$, $\ttt_k$ acts by annihilating all particle configurations whose number of particles is different from $k$.  For particle configurations having $k$ particles, every particle is moved clockwise to the original site of the next particle. Hence,  the particle configuration itself remains fixed by the action of $\ttt_k$ (and it is multiplied with $(-1)^{2(k-1)}q =q$, since a particle has been moved through position $0$).
All the $\ttt_k$ have $\ZZ^N$-degree equal to $(1,\ldots,1)$ and $\ZZ$-degree equal to $N$. Any monomial whose $\ZZ^N$-degree is $(1,\ldots,1)$ occurs as a summand in some central element (after possibly reordering the factors), and the number of summands of $\ttt_k$ equals $\binom{N}{k} =\dim(\bigwedge^k\ground^N)$, see Remark \ref{bem:reconstruct}.
\begin{thm}\label{thm:center}
\begin{enumerate}
\item The $\ttt_k$ are central for all $1\leq k\leq N-1$,  and the center of $\antl{N}$ is generated by $1$ and the $\ttt_k$, $1\leq k\leq N-1$.
\item The subalgebra generated by $\ttt_k$ is isomorphic to the polynomial ring $\ground[q]$  for all $1\leq k\leq N-1$. Moreover $\ttt_k\ttt_\ell=0$ for all $k\neq \ell$. Hence the center of $\antl{N}$ is the subalgebra
$$\zent_N\ =\ \ground\oplus \ttt_1\ground[\ttt_1]\oplus\ldots\oplus \ttt_{N-1}\ground[\ttt_{N-1}]\ \cong\ \frac{\ground[\ttt_1,\ldots,\ttt_{N-1}]}{(\ttt_k\ttt_\ell \mid k\neq \ell)}.$$
\end{enumerate}
\end{thm}
\begin{bew}
\begin{enumerate}
\item The action of $\ttt_k$ on $\V$ is the projection onto the $\antl{N}$-submodule $\ground[q] \otimes  \bigwedge^k\ground^N$ followed by multiplication by $q$. This commutes with the action of every  other element of $\antl{N}$. Since $\V$ is a faithful module, $\ttt_k$ commutes with any element of $\antl{N}$.
As we have seen in Proposition \ref{prop:central}, any central element $c$ without constant term acts on the summand 
$\ground[q] \otimes  \bigwedge^k\ground^N$ via multiplication by some polynomial  $p^c_{k}\in q\ground[q]$. Once again using the faithfulness of $\V$, we get that
$c = \sum\limits_{k=1}^{N-1} p_k^c(\ttt_k)$.
\item 
Recall that $\ground[q] \otimes  \bigwedge^k\ground^N$ is a free $\ground[q]$-module of rank ${N \choose k}$. Since $\ttt_k$ acts by multiplication with $q$ on that module, the subalgebra of $\antl{N}$ generated by $\ttt_k$ must be isomorphic to the polynomial ring $\ground[q]$. Since $\ahj \ahi=0$ for all $\seq{j}\neq \seq{i}$, we get $\ttt_k \ttt_\ell=0$ for $k \neq \ell$, as they consist of pairwise different summands.
\end{enumerate} 
\end{bew}

Theorem \ref{thm:center} enables us to describe the $\ground$-algebra $\End_{\antl{N}}(\W)$  of $\antl{N}$-endomorphisms of the space of nontrivial particle configurations $\W :=\bigoplus\limits_{k=1}^{N-1}\left(\ground[q] \otimes  \bigwedge^k\ground^N\right)   \subset\ \V$. We first observe that on $\W$ multiplication by $q$ is given by the action of a central element in $\zent_N$, therefore it is justified to speak about $\ground[q]$-linearity of a $\antl{N}$-endomorphism of $\W$.
\begin{lemma}\label{lem:zlin}  $\End_{\antl{N}}(\W)\subset\End_{\ground[q]}(\W)$, hence 
any $\antl{N}$-module endomorphism $\varphi$ of $\W$ is $\ground[q]$-linear.   \end{lemma}
\begin{bew}
Observe that $\sum_{k=1}^{N-1}\ttt_k\in\antl{N}$ acts by multiplication by $q$ on every element in $\W$. Therefore multiplication by $q$ commutes with the application of every $\varphi \in \End_{\antl{N}}(\W)$.
\end{bew}
\begin{prop}\label{prop:end}
The endomorphism algebra $\End_{\antl{N}}(\W)$ is isomorphic to a direct sum of $N-1$ polynomial algebras $\ground[T_1]\oplus\ldots\oplus\ground[T_{N-1}]$.
\end{prop}
\begin{bew}
The proof is very similar to the one of Proposition \ref{prop:central}.
First we show that $\varphi(v(\seq{i}))$ is a $\ground[q]$-linear multiple of $v(\seq{i})$ for any $\varphi\in\End_{\antl{N}}(\W)$ and any increasing sequence $\seq{i}$. 
This statement holds if and only if $\pm q\varphi(v(\seq{i}))\in \ground[q]\,v(\seq{i})$. Indeed, by Lemma \ref{lem:tech3} and Lemma \ref{lem:zlin} we get
$$\pm q\varphi(v(\seq{i}))\ =\ \varphi(\pm q v(\seq{i}))\ =\ \varphi(a(\ha{\seq{i}})v(\seq{i}))\ =\ a(\ha{\seq{i}})\varphi( v(\seq{i}))\ \in \ground[q]\, v(\seq{i}).$$
Therefore, we can write $\varphi(v(\seq{i}))=p_{\subsci}\cdot v(\seq{i})$ for some polynomial $p_\subsci \in\ground[q]$. Note that this implies 
$$\End_{\antl{N}}\left(\bigoplus\limits_{k=1}^{N-1}\Big(\ground[q] \otimes \textstyle{\bigwedge^k}\ground^N\Big)\right)\ =\ \bigoplus\limits_{k=1}^{N-1}\left(\End_{\antl{N}}\Big(\ground[q] \otimes \textstyle{\bigwedge^k}\ground^N\Big)\right).$$
What remains is to show that these polynomials only depend on the number of particles in $\seq{i}$, in other words there
exists $p_k \in \ground[q]$ so that  $p_\subsci = p_k$ for all $\seq{i}$ with  $ \vert \seq{i} \vert = k$.  Again it suffices to show this for two sequences $\seq{i}$, $\seq{i}\p$ of length $k$ which differ in exactly one entry. 
So say $i_s=i$, $i_s\p=i+1$, and $i_\ell=i_\ell\p$ for all $\ell\neq s$, for some $1\leq s\leq k$ and $i\in\ZZ/N\ZZ$. 
When $1\leq i\leq N-1$,  
$$
p_{\subsci\p}\,v(\seq{i}\p)\ =\ \varphi(v(\seq{i}\p)) \ =\ \varphi(a_i v(\seq{i}))  \ =\ a_i \varphi( v(\seq{i})) \ = \ a_i (p_\subsci\,v(\seq{i}) ) \ = \ p_\subsci \,v(\seq{i}\p),$$
and when $i = 0$, 
$$(-1)^{k-1}q p_{\subsci\p}\, v(\seq{i}\p)\ =\ (-1)^{k-1}q\varphi( v(\seq{i}\p)) \ =\ \varphi(a_0 v(\seq{i}))  \ =\ a_0 \varphi(v(\seq{i})) \ = \ a_0 (p_{\subsci}\,v(\seq{i})) \ = \ (-1)^{k-1}q p_{\subsci}\,v(\seq{i}\p).$$
Hence we can write $\varphi = \sum_{k=1}^{N-1} p_k \pi_k$ where $\pi_k$ is  the projection onto $\ground[q]\otimes \bigwedge^k\ground^N$, and we get that
$$\End_{\antl{N}}\left(\ground[q] \otimes \textstyle{ \bigwedge^k}\ground^N\right)\ =\ \ground[T_k],$$
where $T_k$ denotes the multiplication action of the central element $\ttt_k$, which is indeed a $\antl{N}$-module endomorphism of $\W$.  Thus, $\End_{\antl{N}}(\W)$ is isomorphic to a direct sum of polynomial algebras as claimed. 
\end{bew}

\begin{bem}
The arguments in the proof of Proposition \ref{prop:end} remain valid even if we specialize the indeterminate $q$ to some element in $\ground\setminus\{0\}$. In this case, we obtain that the summands $\bigwedge^k\ground^N$ are simple modules and $\End_{\antl{N}}\left(\bigoplus\limits_{k=1}^{N-1}\bigwedge^k\ground^N\right)\ \cong\ \ground^{N-1}$.
For $q=0$, the situation is more complicated:  If $q$ is specialized to zero, the generator $a_0$ acts by zero on the module. The action of $\antl{N}$ factorizes over $\fntl{N}$ and the module $\bigwedge^k\ground^N$ is no longer simple. Instead it has a one-dimensional head spanned by the particle configuration $v(1,\ldots,k)$, and any endomorphism is given by choosing an image of this top configuration. It is always possible to map it to itself and to the one-dimensional socle spanned by $v(N-k,\ldots,N)$, but in general there are more endomorphisms. For example, in $\bigwedge^4\ground^8$, the image of $v(1,2,3,4)$ may be any linear combination of $v(1,2,3,4)$, $v(2,3,4,8)$, $v(3,4,7,8)$, $v(4,6,7,8)$ and $v(5,6,7,8)$, so that $\End_{\antl{8}}\left(\bigwedge^4\ground^8\right)$ is  5-dimensional.
\end{bem}

\section{The affine nilTemperley-Lieb algebra is finitely generated over its center}\label{sec:finite}
The affine nilTemperley-Lieb algebra is infinite dimensional when $N \geq 3$; however, the following finiteness result  holds:
\begin{thm}\label{thm:fin}
The algebra $\antl{N}$ is finitely generated over its center.
\end{thm}
\begin{bew}
Given an arbitrary monomial $a(\ul{j})\in\antl{N}$, we first factor it  as $a(\ul{j\p})\cdot a(\ul{j}^{(0)})$ in the following way:  Take the minimal particle configuration $\seq{j} = \{1 \leq j_1 < \ldots < j_k \leq N\}$ on which the monomial $a(\ul{j})$ acts nontrivially. The monomial $a(\ul{j})$ moves all of the particles  by at least one step, because the particle configuration was assumed to be minimal. Using the faithfulness of the representation, we know that we may reorder the monomial $a(\ul{j})$ so that first each particle is moved one step clockwise, and afterwards the remaining particle moves are carried out. Hence,  we may choose some factorization
$a(\ul{j})=a(\ul{j\p})\cdot a(\ul{j}^{(0)})$,  where  $\ul{j}^{(0)}$ is a sequence obtained by permuting $j_1, \dots, j_k$.   The remaining particle moves are carried out by $a(\ul{j\p})$. In Section \ref{sec:appendix},  this decomposition is explicitly constructed (not using the faithful representation).
Next,  we want to find an expression of the form  
$$a(\ul{j})\ =\ a_\text{fin}\cdot \ttt_k^n\cdot a(\ul{j}^{(0)}),$$ 
where $a_\text{fin}$ is a monomial of some subalgebra ${}^i\fntl{N}$ of $\antl{N}$, $\ttt_k^n$ is  in the center of $\antl{N}$, and $a(\ul{j}^{(0)})$ is the above factor. 
Here \begin{equation} {}^i\fntl{N} = \langle a_0,\ldots,a_{i-1},a_{i+1},\ldots,a_{N-1}\rangle\end{equation} is a copy of the finite nilTemperley-Lieb algebra $\fntl{N}$ sitting in $\antl{N}$.
To accomplish this,  we have to subdivide the action of $a(\ul{j})$ on the particle configuration $\seq{j} = \{j_1< \ldots< j_k\}$ one more time. There are two cases:
\begin{enumerate}
\item There is an index $i$ not appearing in $\ul{j\p}$: \ In this case, $a(\ul{j\p})$ is an element of ${}^i\fntl{N}$ and we are done.
\item All indices appear at least $n\geq1$ times in $\ul{j\p}$: \  Let us investigate the action of $a(\ul{j\p})$ on the particle configuration $v(\seq{i})=\ a(\ul{j}^{(0)})\vwj$, where $\seq{i} =\{j_1+1,\ldots,j_k+1\}$. Note that $\seq{i}$  is the minimal particle configuration for $a(\ul{j\p})$.
Each of the particles in $\seq{i}$  is moved by $a(\ul{j\p})$ to the position of the next  particle in the sequence $\seq{i}$,  because there is no index missing (a missing index is equivalent to a particle being stopped before reaching the position of its  successor), before possibly continuing to move along the circle.  Again invoking the faithfulness of the representation, we can rewrite $a(\ul{j\p}) = a(\ul{j\pp})\cdot \ahi^n$,  with the monomial $\ahi$ from Lemma \ref{lem:tech3}. For maximal $n$, the remaining factor $a(\ul{j\pp})$ is an element of ${}^i\fntl{N}$ for some $i$.
Observe that $\ahi^n a(\ul{j}^{(0)})= \ttt_k^n a(\ul{j}^{(0)})$,  which follows immediately from the definition of $\ttt_k$ and Lemma \ref{lem:tech3}. \end{enumerate}
Therefore, we have shown that
$$a(\ul{j})\ =\ a(\ul{j\p})\cdot a(\ul{j}^{(0)})\ =\ a_\text{fin}\cdot \ahi^n\cdot a(\ul{j}^{(0)})\ =\ a_\text{fin}\cdot \ttt_k^n\cdot a(\ul{j}^{(0)}),$$
where $n =0$ in the first case. Since there are only finitely many monomials in $ {}^0\fntl{N},{}^1\fntl{N},\ldots,{}^{N-1}\fntl{N}$ and only finitely many monomials $a(\ul{j}^{(0)})$ such that every index $0,1,\ldots,N-1$ occurs at most once in the sequence $a(\ul{j}^{(0)})$, the affine nilTemperley-Lieb algebra is indeed finitely generated over its center.
\end{bew}

\section{Embeddings of affine nilTemperley-Lieb algebras}\label{sec:inclusion}
In the proof of Theorem \ref{thm:fin},  we have used the $N$ obvious embeddings of $\fntl{N}$ into $\antl{N}$ coming from the $N$ different embeddings of the Coxeter graph $\mathsf{A}_{N-1}$ into $\t{\mathsf{A}}_{N-1}$. Next  we construct $N$ embeddings of $\antl{N}$ into $\antl{N+1}$. They correspond  to the subdivision of an edge of  $\t{\mathsf{A}}_{N-1}$ by inserting a vertex on the edge  to obtain $\t{\mathsf{A}}_N$. 
\begin{thm}\label{thm:emb}
For any number $0\leq m \leq N-1$, there is a unital embedding of algebras $\varepsilon_m:\ \antl{N}\ \rar\ \antl{N+1}$ given by
\begin{equation}\label{eq:emap}
a_i\ \mapsto\ \begin{cases}
a_i\quad&\text{for }\ 0\leq i\leq m-1,\\
a_{m+1} a_{m}\quad&\text{for }\ i=m,\\
a_{i+1}\quad&\text{for }\ m+1\leq i\leq N-1.
\end{cases}
\end{equation} 
\end{thm}
\begin{bem}\label{alghom}  
It is not difficult to see that \eqref{eq:emap} defines an algebra homomorphism $\varepsilon_m$  from 
$\antl{N}$ to $\antl{N+1}$  when $N \geq 3$.   Due to the circular nature of the relations, it suffices to check  this for $\varepsilon_0$.  This  amounts to showing the following, since all the other relations are readily apparent.   To avoid confusion, we indicate generators of $\antl{N+1}$ in these calculations by $\t{a}_i$:  
\begin{align*}&(\t{a}_1\t{a}_0)(\t{a}_1\t{a}_0) = \t{a}_1(\t{a}_0\t{a}_1\t{a}_0) = 0, \quad \t{a}_2(\t{a}_1\t{a}_0)\t{a}_2 = (\t{a}_2\t{a}_1\t{a}_2)\t{a}_0 = 0, \quad \t{a}_N(\t{a}_1\t{a}_0)\t{a}_N  = \t{a}_1(\t{a}_N\t{a}_0\t{a}_N) = 0,\\
&(\t{a}_1\t{a}_0)\t{a}_2(\t{a}_1\t{a}_0) =  (\t{a}_1 \t{a}_2)(\t{a}_0\t{a}_1\t{a}_0) = 0,
\quad (\t{a}_1\t{a}_0)\t{a}_N(\t{a}_1\t{a}_0) =  (\t{a}_1\t{a}_0\t{a}_1)(\t{a}_N\t{a}_0) = 0. 
\end{align*} 
\end{bem}
\begin{bem}\label{bem:busstop}
How should one visualize the action of $\varepsilon_m(\antl{N}) \subset\antl{N+1}$ on the particle configurations on a circle with $N+1$ positions?
Except for $a_m$, all generators of $\antl{N}$ are mapped to corresponding generators of $\antl{N+1}$. They will act as before, by moving a particle one step clockwise around the circle. Since $a_m$ is mapped by $\varepsilon_m$  to the product 
$\t{a}_{m+1}\t{a}_m$ in $\antl{N+1}$, it will move a particle from $m$ to $m+2$, ignoring position $m+1$, as depicted below.
\begin{figure}[H]\label{fig:emb}
\begin{center}
\begin{tikzpicture}[scale=0.6]

\path[use as bounding box] (-3.5,-1.5) rectangle (3.5,2.5);
\draw[thick] (0,0) circle (2cm);

\foreach \x in {225,270,...,360}
\draw[very thick] (\x:1.9cm) -- (\x:2.1cm);


\foreach \x in {45,90,135}
\draw[very thick] (\x:1.9cm) -- (\x:2.1cm);

\foreach \x in {0,...,7}
\draw ({90-\x*45}:1.5cm) node {$\x$};

\foreach \x in {1}
\fill ({90-\x*45}:2.35cm) circle (1.2mm);

\foreach \x in {4}
\draw[thick, densely dotted] ({90-\x*45}:2.35cm) circle (1.2mm);

\foreach \x in {4,7,8}
\draw[thick, ->] ({90-\x*45}:2.7cm) arc ({90-\x*45}:45-\x*45:2.7cm);

\draw[thick, dotted, ->] ({90-5*45}:2.7cm) arc ({90-5*45}:45-6*45:2.7cm);
\draw (70:3cm) node {$\cdot q$};

\end{tikzpicture}
\end{center}
\caption{$\varepsilon_5(\antl{7})\subset\antl{8}$: The action of $\varepsilon_5(a_0a_6a_5a_4)
=\t{a}_0\t{a}_7\t{a}_6\t{a_5}\t{a}_4$ on the particle configuration $v(4)$.}
\end{figure}
\end{bem}

Next we  introduce a basis of $\antl{N}$ that will enable us to see directly that
these homomorphisms are embeddings.  The basis  has a simple description in terms of the graphical representation $\V$ from Section \ref{sec:rep}. 
For any two particle configurations with $1\leq k\leq N-1$ particles  corresponding to the increasing sequences 
$\seq{i}=\{1 \leq  i_1< \ldots < i_k \leq N\}$ and $\seq{j}=\{1 \leq j_1< \ldots< j_k \leq N\}$,  there is a monomial in $\antl{N}$ moving particles at the positions $\seq{j}$ to the positions $\seq{i}$. We require that every particle from 
$\seq{j}$ is moved by at least one step, but we do not prescribe explicitly which of the $j$'s is mapped to which of the $i$'s. For $\seq{i}\neq \seq{j}$,  take $e_{\subsci\subscj}$ to be the monomial such that the power of $q$ in $e_{\subsci\subscj} \vwj =\pm q^\ell \vwi$ is minimal (under the assumption that every particle from $\seq{j}$ must be moved). By faithfulness of the graphical representation, $e_{\subsci\subscj}$ is uniquely determined. For $\seq{i}=\seq{j}$, we have $e_{\subsci\subsci}=\ahi$, the special monomial defined in Section \ref{sec:rep}, hence $e_{\subsci\subsci}\vwi=\pm q\vwi$. Observe that one can write $\ttt_k=\sum_{\mid \subsci \mid=k} e_{\subsci\subsci}$, where the sum runs over all possible increasing sequences  $\seq{i}$ of length $k$, and that $\ttt_k^\ell e _{\subsci \subscj}$ is a monomial,  since all  but one summand vanish for $k=|\seq{i}|$.
The condition that $e_{\subsci \subscj}$ moves all particles from $\seq{j}$ by at least one step guarantees that it acts by zero on all particle configurations with fewer particles than $|\seq{i}|=|\seq{j}|$. 
For example, when $N=7$,
$$e_{(2)(1)} = a_1,\quad e_{(0,2)(0,1)} = a_6a_5a_4a_3a_1a_2a_0a_1.$$ 
(Note that $a_1$ moves $v(0,1)$ to  $v(0,2)$, but this doesn't satisfy the requisite property that 
all the particles must be moved by at least one step.)
If we apply the factorization of monomials from Theorem \ref{thm:fin} to $e_{\subsci \subscj}$, the minimality condition implies that $e_{\subsci \subscj}=a_\text{fin}\cdot 1\cdot a(\ul{j}^{(0)})$,  where if $\seq{j} = \{ j_1 < \ldots < j_k\}$, then 
$\ul{j}^{(0)}$ is a sequence obtained by permuting the elements of $\seq{j}$. 

\begin{thm}\label{thm:easybasis}
The set of monomials
$$\{1\}\ \cup\ \{\ttt_k^\ell  e_{\subsci \subscj} \mid \ell \in\ZZ_{\geq0},\ 1\leq |\seq{i}|=|\seq{j}|=k\leq N-1\}$$
defines a $\ground$-basis of the affine nilTemperley-Lieb algebra $\antl{N}$.
\end{thm}
\begin{bew}
First,  observe that $\ttt_k^\ell e_{\subsci \subscj}$ is indeed a monomial since $|\seq{i}|=k$.
We show that the elements $\ttt_k^\ell e_{\subsci \subscj}$ act $\ground$-linearly independently on the graphical representation $\V=\bigoplus\limits_{k=0}^N \left(\ground[q] \otimes \bigwedge^k \ground^N\right)$. \  By definition, the monomial 
$e_{\subsci \subscj}$ acts by zero on summands $\ground[q] \otimes \bigwedge^{k\p} \ground^N$ for $k\p<|\seq{i}|$. On $\ground[q] \otimes \bigwedge^{|\subsci|} \ground^N$,  the matrix representing the action of $\ttt_k^\ell e_{\subsci \subscj}$ relative to the standard basis has exactly one nonzero entry, and this one distinguishes all monomials with the same minimal number of particles $|\seq{i}|=|\seq{j}|$.
From these two observations, the linear independence follows.
On the other hand, given any nonzero monomial in $\antl{N}$, there exists a minimal particle configuration $\seq{j}$ on which it acts nontrivially. Recording the image particle configuration $\seq{i}$ and the power of $q$, we conclude that there is some $\ell$ so that the element $\ttt_k^\ell  e_{\subsci \subscj}$ acts on $\V$ in the same way as the given monomial does.
Due to the faithfulness of this representation (see Theorem \ref{prop:faith} or Section \ref{sec:appendix}), the proposition follows.
\end{bew} 

In Section \ref{sec:appendix}, a basis is constructed using a different approach (without relying on the faithful representation). Both bases are labelled by pairs of particle configurations (pairs of increasing sequences) together with a natural number $\ell$. Up to an index shift in the output configuration $\seq{i}$ and a shift of the natural number $\ell$, the labelling sets agree,  and both bases actually coincide.
\begin{bew}[Theorem \ref{thm:emb}]
We have already noted in Remark \ref{alghom} that $\varepsilon_m$ is an algebra homomorphism. Using Remark \ref{bem:busstop}, observe that the monomial $e_{\subsci \subscj}\in\antl{N}$ is mapped to a monomial $\tilde e_{\subsci\p \subscj\p}\in\antl{N+1}$ (tilde again indicates in $\antl{N+1}$), where the new index sets are obtained by $i\mapsto i$ for $0\leq i\leq m$ and $i\mapsto{i+1}$ for $m+1\leq i\leq N-1$. 
The injectivity follows since basis elements $\left(\sum_{|\mathrm{\textsc{K}}|=k} e_{\mathrm{\textsc{KK}}} \right)^\ell\cdot e_{\subsci \subscj}$ of $\antl{N}$ are mapped to basis elements $\left(\sum_{\mid \mathrm{\textsc{K}}\p\mid=k} \tilde e_{\mathrm{\textsc{K}}\p \mathrm{\textsc{K}}\p}\right)^\ell\cdot \tilde{e}_{\subsci\p \subscj\p}$ of $\antl{N+1}$.
\end{bew}  
\begin{bem}
It is possible to verify this theorem on generators and relations in the language of Section \ref{sec:appendix} without using the graphical description.
The idea is that from a monomial $e_{\subsci\subscj}$, we can read off the sequences 
\begin{align*}
\seq{j}\ &=\ \{ i \mid  \text{ no }a_{i-1}\text{ to the right of }a_{i} \text{ in the monomial }\ e_{\subsci\subscj} \},\\
\seq{i}\ &=\ \{ i\mid \text{ no }a_{i}\text{ to the left of }a_{i-1}\text{ in the monomial } \ e_{\subsci\subscj}\}.
\end{align*}
Now, using Lemma \ref{lem:nonzero} below, one checks that the image of $e_{\subsci\subscj}$ under $\varepsilon_m$ is a nonzero monomial, which must be equal to the monomial $\tilde{e}_{\subsci\p\subscj\p}$ determined by 
\begin{align*}
\{ i \mid \text{ no }\t{a}_{i-1} \text{ to the right of }\t{a}_{i} \text{ in the monomial }\varepsilon_m(e_{\subsci\subscj}) \} &= \seq{j}\p, \\
\{ i \mid  \text{ no }\t{a}_{i} \text{ to the left of }\t{a}_{i-1} \text{ in the monomial }\varepsilon_m(e_{\subsci\subscj}) \}\ &=\ \seq{i}\p.
\end{align*}
\end{bem}
\begin{bem}
Observe that these embeddings work specifically for the affine \emph{nil}Temperley-Lieb algebras but fail for the ordinary Temperley-Lieb algebras. The relation that fails to hold is the braid relation for Temperley-Lieb algebras, i.e. $a_i a_{i\pm1}a_i = a_i$. Interestingly, the relation $a_i^2=\delta a_i$ is respected for $\delta=1$.
\end{bem}
\newpage
\section{A normal form and the faithfulness of the graphical representation}\label{sec:appendix}
In this section,  we prove Proposition \ref{prop:faith} which we recall here:
\begin{propnn}
For $N\geq3$, $\V$ is a faithful $\antl{N}$-module with respect to the action described in Definition \ref{defi:action}.
\end{propnn}

For the proof, we will explicitly prove  the linear independence of the matrices representing the monomials in $\antl{N}$. We proceed in three steps:  (1) First, we define a normal form for the monomials.  (2) Next,  we find a bijection between the monomials and certain pairs of particle configurations together with a power of $q$. In other words, we find a basis for $\antl{N}$ and describe a labeling set.  (3) The final step is the description of the action of a monomial on $\V$
using its matrix realization.  The matrices representing the monomials have a distinguished nonzero entry that is given in terms of the particle configurations and the power of $q$ from the bijection, and for most matrices, this is the only nonzero entry. From this description it will quickly follow that all these matrices are linearly independent.
\subsection*{Some useful facts}
The following lemma characterises nonzero monomials in $\antl{N}$. They correspond to fully commutative elements in $\atl{N}$, see \cite{green}.
\begin{lemma}\label{lem:nonzero} The monomial
$a(\ul{j})\neq 0$ if and only if for any two neighbouring appearances of $a_i$ in $a(\ul{j})$ there are exactly one $a_{i+1}$ and one $a_{i-1}$ in between, apart from possible factors $a_\ell$ for $\ell \neq i-1,i,i+1$ (indices to be understood modulo $N$). \end{lemma}

According to this result, two consecutive $a_i$ have to enclose $a_{i+1}$ and $a_{i-1}$, i.e. $a_i\ldots a_{i\pm1}\ldots a_{i\mp1}\ldots a_i$, with the dots being possible products of $a_\ell$'s with $\ell \neq i\pm1,i$.
This lemma is a special case of \cite[Lem.~2.6]{green};  here is a quick proof for the convenience of the reader.

\begin{bew} The monomial
$a(\ul{j})$ is zero if and only if we can bring two neighbouring factors $a_i$ together so that we obtain either $a_i^2$ (`square') or $a_ia_{i\pm 1}a_i$ (`braid'). But expressions of the form $a_i\ldots a_{i\pm1}\ldots a_{i\mp1}\ldots a_i$ cannot be resolved this way by commutativity relations. On the other hand, if  there are  two neighbouring factors $a_i$ with either none or only one of the terms $a_{i\pm1}$ in between, we immediately get either $a_i^2$ or $a_ia_{i\pm 1}a_i$. If there are at least two factors $a_{i+1}$ (or $a_{i-1}$) in between the two $a_i$, one can repeat the argument: Either we can create a square or a braid, or we have at least two factors of the same kind in between. In the case of a square or a braid we are done; otherwise we pick two neighbouring $a_{i+k}$ in the $k^{\text{th}}$ step of the argument. Since we always consider the space in between two neighbouring factors $a_i,a_{i+1},\ldots,a_{i+k}$, none of the previous $a_i,a_{i+1},\ldots,a_{i+k-1}$ occur
 s between the two 
neighbouring $a_{i+k}$. Unless we found a square or a braid in an earlier step, we end up in step $N-1$ with a subexpression of the form  $a_{r}a_{r\pm1}^ma_{r}$ which is zero for any $m\geq 0$.
\end{bew}   
 
\begin{defi}\label{def:order}
For any $i\in\{0,1,\ldots,N-1\}$,  we define a (clockwise) order \  $\stackrel{i}{\prec}$ \ on the set $\{0,1,\ldots,N-1\}$ starting at $i$  by 
$$i\,\stackrel{i}{\prec}\,i+1\,\stackrel{i}{\prec}\,  \ldots\,\stackrel{i}{\prec} i +N-1.$$ 
\end{defi}  

\subsection*{Step 1: A normal form}
Given an arbitrary nonzero monomial $a(\ul{j})$ in $\antl{N}$, reorder its factors according to the following algorithm (as usual, the indices are considered modulo $N$):
\begin{enumerate}
\item
Find all factors $a_i$ in $a(\ul{j})$ with no $a_{i-1}$ to their right. We denote them by $a_{i_1},\ldots,a_{i_k}$, ordered according to their appearance in $a(\ul{j})$; in other words, $a(\ul{j})$ is of the form
$$a(\ul{j})\ =\ \ldots a_{i_1}\ldots a_{i_2} \ldots \ \ldots a_{i_k}.$$
\item
Move the $a_{i_1},\ldots,a_{i_k}$ to the far right, without changing their internal order,
$$a(\ul{j})\ =\ a(\ul{j\p})\cdot (a_{i_1}a_{i_2} \ldots a_{i_k})\ =\ a(\ul{j\p})\cdot a(\ul{j}^{(0)})$$
for $\ul{j}^{(0)} = (i_1,\ldots,i_k)$ and some sequence $\ul{j\p} =( \ul{j}\text{ with }i_1,\ldots,i_k\text{ removed})$. This is possible because 
\begin{enumerate}
\item by assumption, there is no $a_{i-1}$ to the right of an $a_i$ in this list; 
\item if  for some $i$,  $a_{i+1}$ occurs to the right of some $a_i$, then  either $a_i\ldots a_{i+1}\ldots a_i$ would occur 
as a subword without $a_{i-1}$ in between, hence $a(\ul{j})=0$, or else $a_{i+1}$ does not have $a_i$ to its right,  so it is one of the $a_{i_1},\ldots,a_{i_k}$ itself, and will be moved to the far right of $a(\ul{j})$, too;
\item  $a_i$ commutes with all  $a_\ell$ for $\ell \neq i-1, i+1$.
\end{enumerate}
\item
Repeat for $a(\ul{j\p})$ until we get 
$$a(\ul{j})\ =\ a(\ul{j}^{(m)})\cdot a(\ul{j}^{(m-1)})\cdot\ldots\cdot a(\ul{j}^{(1)}) \cdot a(\ul{j}^{(0)})$$
for sequences $\ul{j}^{(m)},\ldots,\ul{j}^{(1)}$ obtained successively the same way as described above. Notice:
\begin{itemize}
\item Inside a sequence $\ul{j}^{(n)}$, every index occurs at most once. If two consecutive indices occur within
$\ul{j}^{(n)}$, they are increasingly ordered using the order \ $\stackrel{i_k}{\prec}$ \  from Definition \ref{def:order}.  
\item For two consecutive sequences $\ul{j}^{(n+1)}$, $\ul{j}^{(n)}$ and  for every index $i^{(n+1)}_r$ occurring in $\ul{j}^{(n+1)}$,  we can  find some  index $i^{(n)}_s$ in $\ul{j}^{(n)}$ such that $i^{(n+1)}_r= i^{(n)}_{s}+1$. 
\item From that property, it also follows that the length of $\ul{j}^{(n+1)}$ is less or equal than the length of $\ul{j}^{(n)}$. 
\end{itemize}
\item Reorder the factors $a(\ul{j}^{(m)}),\ldots, a(\ul{j}^{(1)}), a(\ul{j}^{(0)})$ internally:
\begin{enumerate}
\item Start with $a(\ul{j}^{(0)})$. 
There is some $0\leq \hat{\imath}\leq N-1$ which does not occur in $\ul{j}^{(0)}$, but $\hat{\imath}-1$ occurs.
For example, this is satisfied by $\hat{\imath} =i_k+1$, as  $i_k$ occurs in $\ul{j}^{(0)}$ and is to the right of every other factor of $a(\ul{j})$.  Choose the largest such $\hat{\imath}$ (with respect to the usual order).
Then we can move $\hat{\imath}-1$ to the very right of the sequence $\ul{j}^{(0)}$, because  $\hat{\imath}$ is not present, and $\hat{\imath}-2$ may only occur to the left of $\hat{\imath}-1$ due to the construction of $\ul{j}^{(0)}$. We proceed in the same way with those indices $\hat{\imath}-2, \hat{\imath}-3,\ldots, \hat{\imath}-(N-1)$ that appear in $\ul{j}^{(0)}$.
The result is a reordering of the sequence $\ul{j}^{(0)}$ so that it is increasing from left to right with respect to \ $\stackrel{\hat{\imath}}{\prec}$.
\item Repeat with all other factors $a(\ul{j}^{(1)}), a(\ul{j}^{(2)}), \ldots,  a(\ul{j}^{(m)})$ taking  as the initial right-hand index of the sequence  $\hat \imath,\hat\imath+1,\ldots, \hat\imath+m-1$ respectively, and reordering within each 
$a(\ul{j}^{(n)})$ so that the indices are increasing from left to right with respect to \ $\stackrel{\hat{\imath}+n}{\prec}$.
\end{enumerate}
\end{enumerate}

\begin{bsp} As an example for $\antl{7}$, suppose $a(\ul{j})\ =\ a(6\ 4\ 2\ 1\ 3\ 5\ 4\ 2\ 0\ 6\ 1\ 3\ 2\ 5)$. (We omit the commas
to simplify the notation.)\\
\begin{center}
\begin{tabular}{ll}
Find all $a_i$ without $a_{i-1}$ to their right: 
& $a(6\ 4\ 2\ 1\ 3\ 5\ 4\ 2\ 0\ 6\ \uuline{1}\ 3\ \uuline{2}\ \uuline{5})$\\
& \\
Move them to the far right, and & $a(6\ 4\ 2\ 1\ 3\ 5\ 4\ 2\ 0\ 6\ 3)\cdot a(1\ 2\ 5)$\\
don't change their internal order: & \\
\\
Repeat: & $a(6\ 4\ 2\ 3\ 5\ 4\ 1\ \uuline{2}\ 0\ \uuline{6}\ \uuline{3})\cdot a(1\ 2\ 5)$\\
& $a(6\ 4\ 2\ 3\ 5\ 4\ 1\ 0)\cdot a(2\ 6\ 3)\cdot a(1\ 2\ 5)$\\
& $a(6\ 4\ 2\ \uuline{3}\ 5\ \uuline{4}\ 1\ \uuline{0})\cdot a(2\ 6\ 3)\cdot a(1\ 2\ 5)$\\
& $a(6\ 4\ 2\ 5\ 1)\cdot a(3\ 4\ 0)\cdot a(2\ 6\ 3)\cdot a(1\ 2\ 5)$\\
& $a(6\ \uline{4}\ 2\ \uuline{5}\ \uuline{1})\cdot a(3\ 4\ 0)\cdot a(2\ 6\ 3)\cdot a(1\ 2\ 5)$\\
& $a(6\ 2)\cdot a(4\ 5\ 1)\cdot a(3\ 4\ 0)\cdot a(2\ 6\ 3)\cdot a(1\ 2\ 5)$.\\
\\
With the right-hand indices of the  $a(\ul{j}^{(n)})$,  $n \geq 1$,   & $a(6\ 2)\cdot a(4\ 5\ 1)\cdot a(3\ 4\ 0)\cdot a(2\ 3\ 6)\cdot a(1\ 2\ 5)$. \\
arranged according to $\ha \imath+ m-1 \,\stackrel{\ha \imath}{\succ}\, \ldots \,\stackrel{\ha \imath}{\succ}\, \ha \imath +1 \,\stackrel{\ha \imath}{\succ}\, \ha \imath=6$ \\
 from left to right,  reorder the factors in each $a(\ul{j}^{(n)})$ & \\
increasingly with respect to $\stackrel{\ha \imath + n}{\prec}$ from left to right:& 
\end{tabular}
\end{center}
\end{bsp}

As a shorthand notation, in the following we often identify the index sequence $\ul{j}$ with $a(\ul{j})$ (and manipulate $\ul{j}$ according to the same relations as $a(\ul{j})$) as demonstrated in the following example.
\begin{bsp}\label{ex:running}
Let $N=6$.
\begin{align*}
(5\ 1\ 2\ 3\ 0\ 4\ 1\ 5\ 0\ 2\ 3\ 1\ 4\ 5\ 0\ 2\ 3\ 1\ 4\ 2)\ 
&=\ (1)(5\ 0\ 2)(3\ 4\ 5\ 1)(2\ 3\ 4\ 0)(1\ 2\ 3\ 5)(0\ 1\ 2\ 4)\\
&=\ (1\quad 5\ 0\ 2\quad 3\ 4\ 5\ 1\quad 2\ 3\ 4\ 0\quad 1\ 2\ 3\ 5\quad 0\ 1\ 2\ 4).
\end{align*}
\end{bsp}
\begin{lemma}\label{lem:welldef}
Let $a(\ul{j})$ be a nonzero monomial in $\antl{N
}$ with factors indexed by elements in $\ZZ/N\ZZ$. Let $a(\ul{j}^{(m)}),a(\ul{j}^{(m-1)}),\ldots, a(\ul{j}^{(1)}),a(\ul{j}^{(0)})$ be the monomials constructed by the algorithm above. 
\begin{enumerate}
\item
The equality $a(\ul{j})\ =\ a(\ul{j}^{(m)})a(\ul{j}^{(m-1)})\ \cdots \ a(\ul{j}^{(1)})a(\ul{j}^{(0)})$ holds in $\antl{N}$.
\item
Given any two representatives $a(\ul{j})$, $a(\ul{j^\#})$ of the same element in $\antl{N}$, the above algorithm creates the same representative
\ $a(\ul{j}^{(m)}) a(\ul{j}^{(m-1)})\ \cdots \ a(\ul{j}^{(1)}) a(\ul{j}^{(0)})$ for both $a(\ul{j})$ and $a(\ul{j^\#})$.
\end{enumerate}
\end{lemma}
\begin{bew}
\begin{enumerate}
\item
The algorithm never interchanges the order of two factors $a_i$, $a_{i\pm 1}$ with consecutive indices within $a(\ul{j})$. Hence,  the reordering of  the factors of $a(\ul{j})$ uses only the commutativity relation $a_i a_j = a_j a_i$ for $i-j \neq \pm 1\text{ mod }N$ of $\antl{N}$.
\item
Two monomials $a(\ul{j})$, $a(\ul{j^\#})$ in $\antl{N}$ are equal if and only if they only differ by applications of  commutativity relations $a_i a_j = a_j a_i$ for $i-j \neq \pm 1\text{ mod }N$, hence,  if and only if  they contain the same number of factors $a_i$ for each $i$ and the relative position of each $a_i$ and $a_{i\pm 1}$ is the same. Since the outcome of the algorithm depends only on the relative positions of consecutive indices, the resulting decomposition $a(\ul{j}^{(m)})a(\ul{j}^{(m-1)})\ \cdots \ a(\ul{j}^{(1)})a(\ul{j}^{(0)})$ is the same.
\end{enumerate}
\end{bew}

We have shown the following.  In stating this result and subsequently, whenever we refer to monomials in normal form, 
we assume the monomial is nonzero and nonconstant, in particular the sequence $\ul{j}$ is nonempty. 

\begin{thm}\label{thm:basis}
$\{a(\ul{j})  \ \text{in normal form}\}\cup\{1\}$ is a $\ground$-basis of $\antl{N}$.
\end{thm}
\subsection*{Step 2: Labelling of basis elements}\label{sec:psi}
Given $a(\ul{j})\ =\ a(\ul{j}^{(m)})a(\ul{j}^{(m-1)})\ \cdots \ a(\ul{j}^{(1)})a(\ul{j}^{(0)})$ in the normal form, we call $\ul{j}^{(\ell)}$ the $\ell^\text{th}$ {\it block} of $\ul{j}$,  and a string of indices of maximal length of the form $i_s\in\ul{j}^{(0)},\ i_s+1\in\ul{j}^{(1)}, \ i_s+2\in\ul{j}^{(2)},\ldots$ (modulo $N$)  the $s^\text{th}$ {\it strand} of $\ul{j}$.

\begin{bsp}
Let $N=6$, and consider Example \ref{ex:running} once again, where
$$\ul{j}=(1\quad 5\ 0\ 2\quad 3\ 4\ 5\ 1\quad 2\ 3\ 4\ 0\quad 1\ 2\ 3\ 5\quad 0\ 1\ 2\ 4).$$
The blocks are $\ul{j}^{(0)}=(0124)$, $\ul{j}^{(1)}=(1235)$, $\ul{j}^{(2)}=(2340)$, $\ul{j}^{(3)}=(3451)$, $\ul{j}^{(4)}=(502)$, and $\ul{j}^{(5)}=(1)$.
The strands are $[3210]$, $[54321]$, $[105432]$ and $[21054]$. In particular, strands (and blocks) can have different lengths, but the longest strand has length $m=6$.
\end{bsp}

Each monomial $a(\ul{j})\in\antl{N}$ determines two sets $\inj, \outj$ and an integer $\ellj \in\ZZ_{\geq 0}$ as follows:
\begin{align*}
\inj \ &=\ \{i \in\{0,1,\ldots,N-1\}\ |\ \text{no }i-1\text{ to the right of }i\text{ in }\ul{j}\}\\
\outj\ &=\ \{i \in \{0,1,\ldots,N-1\}\ |\ \text{no } i+1\text{ to the left of } i\text{ in }\ul{j}\}\\
\ellj \ &= \text{ the number of zeros  in }\ul{j}. 
\end{align*}
These are well defined because, as in the proof of Lemma \ref{lem:welldef}, any element of $\antl{N}$ is uniquely determined by the number of factors $a_i$ and the relative position of each $a_i$ and $a_{i\pm 1}$, for all $i$. 
The set  $\seq{i}^{\text{in}}_{\ul{j}}$ equals the underlying set of $\ul{j}^{(0)}$ in the normal form from the algorithm above. All strands of $\ul{j}$  begin with an element in $\inj$ and end with an element from $\outj$.   

The goal of this subsection is to show  
\begin{prop}\label{prop:psi}
The mapping 
\begin{align}\label{eq:psidef}
\psi:\ \{a(\ul{j})\in\antl{N}\text{  in normal form}\}\ &\rar\ {\cal P}_N \times {\cal P}_N \times \ZZ_{\geq 0}\\
a(\ul{j})\ &\mapsto\ (\inj,\outj,\ellj),  \nonumber 
\end{align}
is injective,  where ${\cal P}_N$ is the power set of $\{0,1,\ldots,N-1\}$. \end{prop}
\begin{bem}
The map $\psi$ is defined so that in the graphical description of the representation $\V$ of $\antl{N}$, the set 
$\inj$ equals the set of positions where $a(\ul{j})$ expects a particle to be. The set $\outj$ equals  the set of positions where $a(\ul{j})$ moves the particles from $\inj$,  but each one is translated by 1, that is, 
$$a(\ul{j})\text{ applied to a particle at }i\in \inj  \text{ gives a particle at }j+1\text{ for some } j\in \outj.$$
The map $\psi$ is far from being surjective. An obvious constraint is that $\vert \inj \vert = \vert \outj \vert$, and furthermore, for some pairs $(\inj,\outj)$, one can only obtain sufficiently large values $\ellj$.
\end{bem}

To ease the presentation, we start by proving injectivity of the restriction $\psi_0$ of $\psi$
to those monomials $a(\ul{j})$ in normal form whose first element $i_1$ of $\ul{j}^{(0)}$ is $0$. The proof itself  will amount to counting indices.
\begin{prop}\label{prop:psi0}
$$\psi_0:\ \{a(\ul{j})\in\antl{N}\text{  in normal form, with }i_1=0\}\ \rar\  {\cal P}_N \times {\cal P}_N \times \ZZ_{\geq 0}, 
\qquad a(\ul{j})\mapsto (\inj,\outj,\ellj) $$ 
is an injective map.
\end{prop}

Before beginning the proof of this result, we note that for monomials $a(\ul{j})$ with $i_1=0$, the inequality  $i_k<N-1$ must hold
in $\inj$,  since $i_1 = 0$ implies that $i_1-1=N-1$ is not an element of $\inj$.
 Consequently, the  ordering of the indices in $\inj$ agrees with the natural ordering of $\ZZ$, 
so we can regard  $(\inj, <)$ as a subset of $(\ZZ,<)$ and replace the modular index sequence $\ul{j}$ by an integral index sequence $\ul{j}^\ZZ$ such that $\ul{j}^\ZZ (\text{ mod } N) = \ul{j}$.
\begin{defi}\label{defi:int}
Assume  $\ul{j}=\ul{j}^{(m)}\cdot\ldots\cdot  \ul{j}^{(1)} \cdot \ul{j}^{(0)}$ is a normal form sequence with $\ul{j}^{(0)}= \{0=i_1<\ldots < i_k < N-1\}$ and $\ul{j}^{(n)}=(i_{h_1}+n,\ldots,i_{h_{k(n)}}+n)\subseteq (i_1+n,\ldots,i_k+n)$,  where indices in  $\ul{j}^{(n)}$ are modulo $N$ and $1\leq k(n)\leq k$ for all $1\leq n\leq m$.  The {\it integral normal form sequence for $\ul{j}$} is 
$$\ul{j}^\ZZ=(\ul{j}^{(m)})^\ZZ\cdot\ldots\cdot (\ul{j}^{(1)})^\ZZ \cdot \ul{j}^{(0)}\quad\text{where}\quad (\ul{j}^{(n)})^\ZZ\ :=\ (i_{h_1}+n,\ldots,i_{h_{k(n)}}+n)\in\ZZ^{k(n)}$$
for $n = 1,\dots,m$. 
\end{defi}

\begin{bsp}
We continue Example \ref{ex:running} with $N=6$.   
\begin{align*} \text{If} \ \  \ul{j} &=(1\quad 5\ 0\ 2\quad 3\ 4\ 5\ 1\quad 2\ 3\ 4\ 0\quad 1\ 2\ 3\ 5\quad 0\ 1\ 2\ 4),  \\
\text{then}\ \ \ul{j}^\ZZ&=(7\quad 5\ 6\ 8\quad 3\ 4\ 5\ 7\quad 2\ 3\ 4\ 6\quad 1\ 2\ 3\ 5\quad 0\ 1\ 2\ 4).
\end{align*}
\end{bsp}

Our proof of Proposition  \ref{prop:psi0} will hinge upon the following technical lemma.   

\begin{lemma}\label{lem:ineq} Let $\ul{j}^\ZZ$ be the integral normal form sequence for $\ul{j}$ and let   
$[i_s,\ldots, i_s+n_s]$ for $s=1,\dots,k$ be the strands of  $\ul{j}^\ZZ$.   Assume $i_1 = 0$.    Then 
\begin{itemize}
\item[{(a)}]  
$n_1 =i_1+n_1\ <\ i_2+n_2\ <\ \ldots\ <\ i_k+n_k;$ 
\item[{(b)}]  $i_k+n_k< i_1+n_1+N = n_1+N$.
\end{itemize}
\end{lemma}
 
We postpone the proof of this result and proceed directly to proving the proposition.   

\begin{bew}[Proposition \ref{prop:psi0}]  Since the sequence $\ul{j}$ will be fixed throughout the proof, we
will drop the subscript $\ul{j}$ on $\inj \outj, \ellj$.   
To show the injectivity of $\psi_0$, we consider the factorization $\psi_0=\gamma\circ\beta\circ\alpha$ given by 
$$\psi_0:\ a(\ul{j})\ \stackrel{\alpha}{\longmapsto}\ a(\ul{j}^\ZZ)\ \stackrel{\beta}{\longmapsto}\ (\inz,\outz)\ \stackrel{\gamma}{\longmapsto}\ (\seq{i}^\text{in},\seq{i}^\text{out},\ell),$$
where $\inz=\seq{i}^{\text{in}}$ and $\outz=\{i\in\ul{j}^\ZZ\ |\ \text{no } i+1\text{ to the left of }i\}$ similar to the definition of $\seq{i}^\text{out}$. 
The map $\alpha$ replaces indices in $\ZZ/N\ZZ$ by indices in $\ZZ$ as in Definition \ref{defi:int} above. The map $\beta$ is given by  reading off $\outz$ and $\inz$ from $\ul{j}^\ZZ$.
The map $\gamma$  sends  the pair $(\inz,  \outz)$ to a triple 
consisting of the  respective images $\seq{i}^\text{in}, \seq{i}^\text{out}$ modulo $N$ of the pair and the integer $\ell=1+\sum \ell_r$ where $\ell_r=\lfloor\frac{j_r}{N}\rfloor$ for each $j_r\in \outz$. The summand $1$ corresponds to $0=i_1$; all other occurrences of $0$ are counted by $\sum \ell_r$.   Now we check injectivity.
 
The map $\alpha$ is clearly injective since $\ul{j}^\ZZ\ \mapsto\ \ul{j}^\ZZ (\text{ mod } N)$ is a left inverse map.
 
To see that  $\beta$ is injective,  we need to know that $\ul{j}^\ZZ$ can be uniquely reconstructed from $(\inz, \outz)$.
Observe that $\ul{j}^\ZZ$ is determined by knowing all the `strands' $i_s, i_s+1, i_s+2,\ldots, i_s+n_s$ for $1\leq s\leq k$, hence by assigning an element $i_s+n_s\in \outz$ to each $i_s\in \inz$. But  
it follows from Lemma \ref{lem:ineq}\,(a)  that  $i_1+n_1$ must be the smallest element of $\outz$, $i_2+n_2$ the second smallest, etc., so that the element  $i_s+n_s$ is assigned to the $s$th element in $\seq{i}^\text{in}$,
that is,  to $i_s$.    
 
Now to see that  $\gamma$ is injective, we need to recover $(\inz,\outz)$ in a unique way from $(\seq{i}^\text{in},\seq{i}^\text{out},\ell)$.   Write  $\seq{i}^{\text{in}}=\{0=i_1<\ldots<i_k<N-1\}$,  and set $\inz:=\seq{i}^{\text{in}}$. 
By Lemma \ref{lem:ineq}\,(a),  we know that $\outz$ is of the form $(i_1+n_1<\ldots<i_k+n_k)$, and since the elements of $\seq{i}^\text{out}$ have to be equal to the elements of $\outz$  modulo $N$, we can write $i_r+n_r=N\ell_r+d_r$ for $\ell_r=\lfloor\frac{i_r+n_r}{N}\rfloor$ and some $d_r\in \seq{i}^\text{out}$.
Comparing $\ell_r$ and $\ell_s$ for $r<s$, we have 
$$N\ell_r\ \leq\ N \ell_r+d_r\ =\ i_r+n_r\ <\ i_s+n_s\ =\ N\ell_s+d_s\ \leq\ N(\ell_s+1).$$
So $\ell_r< \ell_s+1$, i.e. $\ell _r\leq \ell_s$. Similarly,  we obtain from (b) of Lemma \ref{lem:ineq} that $\ell_k\leq \ell_1+1$.   

As a result, 
$$N \ell_k\ \leq\ N\, \ell_k+d_k \ =\ i_k+n_k\ <\ i_1+n_1+N\ =\ N(\ell_1+1)+d_1\ \leq\ N(\ell_1+2),$$
i.e. $\ell_k<\ell_1+2$. Together we have $\ell_1=\ldots=\ell_s <\ell_{s+1}=\ldots=\ell_1+1$ for some $1<s\leq k$ (where we treat the case $s=k$ by $\ell_1=\ldots=\ell_k$).
Set $\t{\ell}:=\ell_1$. Then
\begin{align*}
i_r + n_r\ &=\ N\,\t{\ell}+d_r\quad&\text{for }1\leq r\leq s,\\
i_r + n_r\ &=\ N(\t{\ell}+1)+d_s \quad&\text{for }s+1\leq r\leq k.
\end{align*}
As a first consequence, 
$$\ell\ =\ 1+\sum_r \ell_r\ =\ 1+k \t{\ell}+(k-s),$$
which determines $\t{\ell}= \lfloor\frac{\ell-1}{k}\rfloor$,  and hence all $\ell_r$,  as well as the index $s$.
Using Lemma \ref{lem:ineq}, we determine that 
$$i_{s+1}+n_{s+1}\ <\ \ldots\ <\ i_k+n_k\ <\ i_1+n_1+N\ <\ \ldots\ <\ i_s+n_s+N,$$
and so
$$N\,(\t{\ell}+1)+d_{s+1}\ <\ \ldots\ <\ N\,(\t{\ell}+1)+d_k\ <\ N\,(\t{\ell}+1)+d_1\ <\ \ldots\ <\ N\,(\t{\ell}+1)+d_s.$$
Therefore, $d_{s+1}<\ldots<d_k<d_1<\ldots<d_s$, which fixes the choice of $d_r$ for all $r$.
We conclude that given $(\seq{i}^\text{in},\seq{i}^\text{out},\ell)$, we can reconstruct $\outz$ by setting $i_r+n_r:=N\, \ell_r+d_r$.
This completes the proof of Proposition  \ref{prop:psi0}.
\end{bew}  

\begin{bew}[Lemma \ref{lem:ineq}]
(a)  Let $\ul{j}^\ZZ$ be a nonempty integral normal form sequence with $0=i_1<\ldots<i_k\leq N-1$ and strands $[i_r,\ldots,i_r+n_r]$ for $1\leq r\leq k$.
Assume that there is some index $1\leq t\leq k-1$ such that 
$i_t+n_t\ \geq\ i_{t+1}+n_{t+1}$.
Since $i_t<i_{t+1}$, we have $n_t> n_{t+1}$. So 
$$\ul{j}^\ZZ\ =\ \ldots \underbrace{(\ldots\ i_{t}+n_{t}\ \ldots)}_{\text{the } n_{t}\text{th bracket}}\ \ldots\ \underbrace{(\ldots\ i_t+n_{t+1}\quad i_{t+1}+n_{t+1}\ \ldots)}_{\text{the }n_{t+1}\text{th bracket}}\ \ldots\ .$$
From $i_t +n_{t+1}\ <\ i_{t+1}+n_{t+1}\ <\ i_t+n_t$ it follows that there is some integer  $n_{t+1}<p\leq n_t$ such that $i_{t+1}+n_{t+1}=i_t+p$ appears in the strand $[i_t,\ldots,i_t+n_t]$, i.e.
$$\ul{j}^\ZZ\ =\ \ldots \underbrace{(\ldots\ i_{t}+n_{t}\ \ldots)}_{\text{the } n_{t}\text{th bracket}}\ \ \ldots \underbrace{(\ldots\ i_{t}+p\ \ldots)}_{\text{the } p \text{th bracket}}\ \ldots\ \underbrace{(\ldots\ i_t+n_{t+1}\quad i_{t+1}+n_{t+1}\ \ldots)}_{\text{the } n_{t+1}\text{th bracket}}\ \ldots\ $$
with $i_t+p=i_{t+1}+n_{t+1}$. But by the definition of the strands,  there is no $i_{t+1}+n_{t+1}+1$ appearing to the left of $i_{t+1}+n_{t+1}$.  Due to Lemma \ref{lem:nonzero},  we know that (even modulo $N$) there is no repetition of $i_{t+1}+n_{t+1}$ to the left. Thus $i_t+p=i_{t+1}+n_{t+1}$ is not possible, and we obtain $i_1+n_1\ <\ i_2+n_2\ <\ \ldots\ <\ i_k+n_k$. 

For (b) of Lemma \ref{lem:ineq},  assume $i_k+n_k\geq i_1+n_1+N$.  It is true generally that $N>i_k$, so we get
$i_k+n_k \geq i_1+n_1+N>i_k+n_1$. Hence $i_1+n_1+N=i_k+b$ for some $n_1<b\leq n_k$, i.e. $i_1+n_1+N$ appears in the strand $[i_k,\ldots, i_k+n_k]$ and we have
$$\ul{j}^\ZZ\ =\ \ldots \underbrace{(\ldots\ i_k+n_k)}_{\text{the } n_k\text{th bracket}}\ \ldots\ \underbrace{(\ldots\ i_k+b\ \ldots)}_{\text{the }b\text{th bracket}}\ \ldots\ \underbrace{(i_1+n_1\ \ldots\ i_k+n_1)}_{\text{the } n_1\text{th bracket}}\ \ldots\ .$$
Here it may be that the $n_k$th bracket and the $b$th bracket coincide, but in any case,  we find that $i_k+b\ =\ i_1+n_1+N\ = \ i_1+n_1\text{ mod } N$, and so $i_k+b$  appears to the left of $i_1+n_1$. By the definition of the strands, there is no $i_1+n_1+1$ to the left of $i_1+n_1$, and from Lemma \ref{lem:nonzero} we deduce that in $\ul{j}=\ul{j}^\ZZ\text{ mod }N$ there is no $i_1+n_1\text{ mod }N$ to the left of $i_1+n_1$ allowed, which leads to a contradiction. Hence $i_k+n_k < i_1+n_1+N$ must hold.  
\end{bew}

Having established  that $\psi$ is injective when restricted to sequences with $i_1 = 0$, we now show the injectivity of $\psi$
in general. 

\begin{bew}[Proposition \ref{prop:psi}]
We have the following disjoint decompositions according to the smallest value $i_1$ 
in  $\ul{j}^{(0)}$ for $\ul{j}$: 
\begin{align*}
\{a(\ul{j}) \ \text{in normal form}\}\ &=\ \coprod_i \{a(\ul{j})\ \text{in normal form, with} \  i_1=i\, \}\\
\{(\inj,\outj,\ellj)\}\ &=\ \coprod_i \{(\inj,\outj,\ellj) \mid 
i_1 = i \in \inj\}\\
\psi\ &=\ \coprod_{i}\left(\psi_i:\ \{a(\ul{j})\ \text{in normal form, with} \  i_1=i\, \}\ \rar\ \{(\inj,\outj,\ellj)
\mid i_1 = i \in \inj  \}\right).
\end{align*}
By Proposition \ref{prop:psi0}, the map $\psi_0:\ a(\ul{j})\mapsto (\inj,\outj,\ellj)$ restricted to those $a(\ul{j})$ with $i_1=0$ is injective.  We argue next  that by an index shift this result is true for all other $\psi_i$. 
Now it follows from Proposition \ref{prop:psi0} that the map  
$$\h\psi_0:\ \{a(\ul{j})\in\antl{N}\text{  in normal form, with }i_1=0\}\ \rar\ \{(\inj,\outj,\h{\ellj})\mid i_1=0 \in
\seq{i}_\text{in}\}$$  
is injective, where $\h{\ellj}$ counts the occurences of $N-i$ in $\ul{j}$. Recall that 
$$\ellj=\sum_r \ell_r+1\ \ \text{and}\ \ \ell_r \text{ is the number of 0 in the }r\text{th strand }\  [i_r,\ldots,i_r+n_r] \ \text{of } \ul{j}\  \text{mod }N.$$
Now observe that we can obtain $\ellj$ from $\h{\ellj}$ as
$$\ellj=\h{\ellj}-\big | \{d_r \in \outj \mid d_r \geq N-i \}\big |+\big | \{ i_r\in\inj \mid i_r>N-i\}\big |+1,$$
which follows from a computation using $\h{\ellj}=\sum_r \h{\ell_r}$ and
\begin{align*}
\h{\ell_r}\ & = \text{ the number of \ $N-i$ \  in the }r\text{th strand } \ [i_r,\ldots,i_r+n_r]\text{ mod }N \\
&=\ \begin{cases}\lfloor\frac{i_r+n_r+i}{N}\rfloor\quad&\text{if }i_r\leq N-i\\
\lfloor\frac{i_r+n_r+i}{N}\rfloor-1\quad&\text{if }i_r> N-i\end{cases}\\
&=\ \begin{cases}\lfloor\frac{N \ell_r+d_r+i}{N}\rfloor\quad&\text{if }i_r\leq N-i\\
\lfloor\frac{N \ell_r+d_r+i}{N}\rfloor-1\quad&\text{if }i_r> N-i\end{cases}\\
&=\ \begin{cases}\ell_r+1\quad&\text{if }i_r\leq N-i\text{ and } d_r+i\geq N\\
\ell_r\quad&\text{if }i_r\leq N-i\text{ and }d_r +i< N\\
\ell_r\quad&\text{if }i_r> N-i\text{ and }d_r+i\geq N\\
\ell_r-1\quad&\text{if }i_r> N-i\text{ and }d_r+i< N.\end{cases} 
\end{align*}

We obtain $\psi_i$ by first shifting the indices of $\ul{j}$ by subtracting $i$ from each index,  $\ul{j}-(i,\ldots,i)$, then applying $\h \psi_0$, and finally shifting the indices from $\inj$ and $\outj$ by adding $i$ to each.   Hence, $\psi_i$ is injective for each $i$, and $\psi$ is injective because the unions are disjoint.
\end{bew}

\subsection*{Step 3: Description and linear independence of the matrices}

Recall that the standard $\ground$-basis of the representation $\V=\bigoplus\limits_{k=0}^N\left(
\ground[q] \otimes \bigwedge^k\ground^N\right)$ is given by  
$$\{ q^\ell \cdot v_{i_1}\wedge\ldots\wedge v_{i_k}\mid \ell \in\ZZ_{\geq0},\ 1\leq i_1<\ldots <i_k\leq N\}$$
where $(i_1,\ldots,i_k)$ is identified with the particle configuration having particles in those positions in the graphical description.
Now we describe with respect to this basis the matrix representing a nonzero monomial $a(\ul{j})\in\antl{N}$ as a $2^N\times2^N$-matrix with entries in $\ground[q]$.
Since $\V$ decomposes as a $\antl{N}$-module into submodules $\ground[q] \otimes \bigwedge^k\ground^N$
for $k = 0,1,\dots, N$, the matrix of $a(\ul{j})$ is block diagonal with $N+1$ blocks $A_0,A_1,\ldots,A_N$, where $A_0=A_N=(0)$ corresponding  to the trivial representation. 
$$
a(\ul{j})\ =\ \left(\begin{array}{ccccc}
0 & 0 &  & \cdots & 0 \\
0 & \begin{array}{|c|}\hline A_1\\ \hline \end{array} &  &  & \vdots\\
 &  & \ddots &  & \\
\vdots &  &  & \begin{array}{|c|}\hline A_{N-1}\\ \hline \end{array} & 0\\
0 & \cdots &  & 0 & 0\\
\end{array}\right)
$$
The block $A_k$ is a $\binom{N}{k}\times\binom{N}{k}$-matrix, with entries from $\ground[q]$ indexed by all possible particle configurations whose number of particles equal to $k$. 

Now fix a nonzero monomial $a(\ul{j})$ in normal form that is specified by the triple $(\inj,\outj,\ellj)$ defined in Step 2. Let $k=\vert \inj\vert$.
All blocks $A_1,\ldots,A_{k-1}$ are zero since $a(\ul{j})$ expects at least $k$ particles. For $r>k$, there might be nonzero blocks (unless the particles from $\inj$ are moved around the whole circle with no position left out, in which case there are no surplus particles allowed. This occurs if $a(\ul{j})$ contains at least every other generator $a_i,a_{i+2},\ldots$).
More importantly, the block $A_k$ has precisely one nonzero entry, and this is given by 
$$ (A_k)_{\subsci^\text{in}_{\ul{j}},\subsci^\text{out}_{\ul{j}}}\ =\ \pm q^{\ellj}.$$
From this we see first that all matrices representing monomials $a(\ul{j})$ in normal form with $\vert \inj\vert=N-1$ are $\ground$-linearly independent: They have only one nonzero entry which is equal to $\pm q^{\ellj}$ at position $(\inj,\outj)$. Furthermore, if all matrices representing monomials $a(\ul{j})$ in normal form with $\vert \inj\vert \geq k$ are $\ground$-linearly independent, then also all matrices representing monomials $a(\ul{j})$ in normal form with $\vert \inj\vert  \geq k-1$ are $\ground$-linearly independent. This follows because the additional monomials $a(\ul{j})$ with $\vert \inj\vert =k-1$ have nonzero entries $(A_{k-1})_{\subsci^\text{in}_{\ul{j}},\subsci^\text{out}_{\ul{j}}}=\pm q^{\ellj}$
in the $(k-1)$th block which is zero for all $a(\ul{j})$ with $\vert \inj\vert \geq k$.
So by induction,  all matrices representing monomials $a(\ul{j})$ in normal form are $\ground$-linearly independent. 
Since all of them have a zero entry in the upper left (and lower right) corner, we may add the identity matrix to the linearly independent set of matrices, and it remains linearly independent. So the representation of $\antl{N}$ on $\V$ is faithful, because according to Theorem \ref{thm:basis}, $\{a(\ul{j})\ \text{in normal form}\}\cup\{1\}$ is a $\ground$-basis of $\antl{N}$.

Section 8 has given a normal form for each monomial and has provided an alternate proof of the faithfulness of the representation of $\antl{N}$ by elementary arguments.

\bibliographystyle{amsalpha}
\addcontentsline{toc}{section}{Literature}
\def\cprime{$'$}
\providecommand{\bysame}{\leavevmode\hbox to3em{\hrulefill}\thinspace}
\providecommand{\MR}{\relax\ifhmode\unskip\space\fi MR }
\providecommand{\MRhref}[2]{%
  \href{http://www.ams.org/mathscinet-getitem?mr=#1}{#2}
}
\providecommand{\href}[2]{#2}

\medskip

 \noindent \textit{\small Department of Mathematics, University of Wisconsin-Madison, Madison, WI 53706, USA}\\
{\small benkart@math.wisc.edu}
\medskip

\noindent \textit{\small   Mathematical Institute, 
University of Bonn, 53115  Bonn, Germany}\\
{\small joanna@math.uni-bonn.de}

\end{document}